\documentclass{article}
\usepackage{fullpage}

\usepackage{amsmath}
\usepackage{amsfonts}
\usepackage{amsthm}
\usepackage{thmtools}
\usepackage{thm-restate}
\usepackage{tikz-cd}

\usepackage{hyperref}
\usepackage[noabbrev]{cleveref}

\newtheorem{lemma}{Lemma}[section]

\newtheorem{proposition}[lemma]{Proposition}
\newtheorem{corollary}[lemma]{Corollary}

\theoremstyle{definition}
\newtheorem{definition}[lemma]{Definition}
\theoremstyle{remark}
\newtheorem{remark}[lemma]{Remark}
\newtheorem{family}{Example Family}

\numberwithin{equation}{section}

\usepackage{macros}
\usepackage{graphicx}

\crefname{section}{Section}{Sections}
\crefname{lemma}{Lemma}{Lemmas}
\crefname{proposition}{Proposition}{Propositions}
\crefname{theorem}{Theorem}{Theorems}
\crefname{corollary}{Corollary}{Corollaries}
\crefname{conjecture}{Conjecture}{Conjectures}
\crefname{definition}{Definition}{Definitions}
\crefname{table}{Table}{Tables}
\crefname{equation}{Equation}{Equations}
\crefname{figure}{Figure}{Figures}
\crefname{family}{Example Family}{Example Families}

\begin{document}

\title{Stable geodesic nets in convex hypersurfaces}
\author{Herng Yi Cheng\\University of Toronto\\\texttt{herngyi.cheng@mail.utoronto.ca}}
\date{}
\maketitle

\begin{abstract}
    We construct convex bodies that can be ``captured by nets.'' More precisely, for each dimension $n \geq 2$, we construct a family of Riemannian $n$-spheres, each with a stable geodesic net, which is a stable 1-dimensional integral varifold. Small perturbations of a stable geodesic net must lengthen it. These stable geodesic nets are composed of multiple geodesic loops based at the same point, and also do not contain any closed geodesic. All of these Riemannian $n$-spheres are isometric to convex hypersurfaces of $\R^{n+1}$ with positive sectional curvature.
\end{abstract}


\section{Introduction}

This paper will present constructions of convex bodies that can be ``captured by nets.'' Formally speaking, these ``nets'' are \emph{stable geodesic nets}. To define them, we begin with the more general notion of a \emph{stationary geodesic net} in a Riemannian manifold $M$, which is an immersion $G : \Gamma \to M$ of a graph $\Gamma$ whose edges are mapped to constant-speed geodesics, such that the outgoing tangent vectors at each vertex sum to zero. (In this context, an immersion of a graph is a continuous map that is a piecewise $C^\infty$ immersion when restricted to each edge.) For example, on a round 2-sphere, the union of three lines of longitude at angles $2\pi/3$ to each other is a stationary geodesic net. 

Among all possible immersions of $\Gamma$ in $M$, stationary geodesic nets are the critical points of the length functional. In geometric measure theory, stationary geodesic nets arise as stationary 1-dimensional integral varifolds \cite{Pitts_StationaryVarifolds,AllardAlmgren_StationaryVarifolds}. They are intimately connected to the study of closed geodesics. Attempts to find closed geodesics using min-max methods, such as searching for critical points of the length functional in the space of one-dimensional flat cycles, may fail as they may yield stationary geodesic nets that do not contain closed geodesics \cite{NabutovskyRotman_GeodesicNets}.

There are very few existence results on stationary geodesic nets. ``Trivial'' stationary geodesic nets can be formed as the union of some closed geodesics. However, most of the existence results to date can guarantee the existence of a stationary geodesic net without being able to tell whether it contains a closed geodesic or not. This is the case for the proof by A.~Nabutovsky and R.~Rotman of the existence of short stationary geodesic nets on any closed manifold \cite{NabutovskyRotman_GeodesicNets,Rotman_Flowers} and the recent proof by Y.~Liokumovich and B.~Staffa that the union of stationary geodesic nets is dense in generic closed Riemannian manifolds \cite{LiokumovichStaffa_GenericDensity}. On the other hand, J.~Hass and F.~Morgan proved that any metric on $\Sp[2]$ sufficiently close to the round metric in the $C^2$ topology has a stationary geodesic net homeomorphic to the $\theta$ graph, which cannot contain any closed geodesic \cite{HassMorgan_GeodesicNetS2}.\footnote{The $\theta$ graph is the graph with two vertices that are connected by three edges.} We are not aware of any proof that there is some closed manifold $M$ of dimension 3 or more such that for all Riemannian metrics $g$ in some open set in the $C^k$ topology for some fixed $2 \leq k \leq \infty$, $(M,g)$ admits a stationary geodesic net which contains no closed geodesics. In this paper we prove the existence of such an open set, in the $C^\infty$ topology, of metrics on $\Sp[n]$ for all $n \geq 3$.

We will study \emph{stable geodesic nets}; they are stationary geodesic nets that are also local minima of the length functional among all immersions of the same graph. However we will work with a stronger notion, to be formally defined in \cref{sec:Definitions}, that corresponds to the notion of non-degenerate critical points in Morse theory. One application of stable geodesic nets lies in a certain approach to find a short closed geodesic in a closed manifold: one can pick a triangulation of the manifold with bounded edge lengths, perform a length-shortening process on the 1-skeleton of each simplex, and show using topological methods that one of of the 1-skeletons must converge to a stable geodesic net that is not a point. However, this approach cannot tell \emph{a priori} whether this stable geodesic net will contain a closed geodesic or not. For this reason, proving the non-existence of certain classes of stable geodesic nets can help to prove length bounds on the shortest closed geodesic in a manifold. For example, I.~Adelstein and F.\,V.~Pallete showed that on positively curved Riemannian 2-spheres, $\theta$ graphs are never stable. They used this result to improve the best-known bound on the length of the shortest closed geodesic on a 2-sphere in terms of the diameter, under the assumption of non-negative curvature \cite{AdelsteinPallete}.

Indeed, sufficiently pinched positive curvature in a manifold can obstruct the existence of stable geodesic nets. A well-known conjecture by H.\,B.~Lawson~Jr. and J.~Simons asserts that 1/4-pinched manifolds that are compact and simply-connected cannot contain stable geodesic nets and other stable varifolds or stable submanifolds \cite{LawsonSimons}. This conjecture has been partially confirmed for certain classes of manifolds satisfying various curvature pinching conditions \cite{ShenXu_QuarterPinchedHypersurfaces,Howard_SuffPinched, HuWei_FifthPinched}. Furthermore, a classical theorem of J.\,L.~Synge forbids even-dimensional, compact and orientable manifolds from having any stable closed geodesics if they have positive sectional curvature---regardless of the degree of curvature pinching \cite{Synge_GeodesicPosCurv}.

We construct examples of stable geodesic nets in positively-curved Riemannian spheres of every dimension; in particular they are immersions of graphs that will consist of multiple geodesic loops based at the same point. We will call them \emph{geodesic bouquets}.

\begin{restatable}[Main result]{theorem}{ThmStableNLoopsPosCurv}
    \label{thm:StableNLoopsPosCurv}
    For every integer $n \geq 3$, there exists a convex hypersurface $M_n$ of $\R^{n+1}$ with positive sectional curvature, such that it contains a stable geodesic bouquet with $n$ loops. There also exists a convex surface $M_2$ in $\R^3$ with positive sectional curvature that contains a stable geodesic bouquet with 3 loops.
    
    Furthermore, for each $n \geq 2$, every Riemannian $n$-sphere whose metric is sufficiently close to that of $M_n$ in the $C^\infty$ topology also contains a stable geodesic bouquet. All of the stable geodesic bouquets produced in this theorem contain no closed geodesic.
\end{restatable}

Our main result contrasts with the theorem of Synge: even-dimensional, compact and orientable manifolds cannot contain stable closed geodesics, but they may contain stable geodesic nets. We also note that this main result disproves a conjecture of R.~Howard and S.\,W.~Wei in every dimension \cite[Conjecture~B]{HowardWei_Hypersurfaces}, which implies in particular that simply-connected closed manifolds with positive sectional curvature cannot contain stable closed geodesics or stable geodesic nets.\footnote{We also highlight a counterexample by W.~Ziller, of a closed homogeneous 3-manifold with positive sectional curvature and a stable closed geodesic \cite[Example~1]{Ziller_HomogSpaceStableClosedGeodesic}.}

All of the stable geodesic bouquets that we construct will be ``elementary'' in the sense that they are essentially constructed using straight line segments in convex polytopes. As a result, most of the constructions can be understood using only Euclidean geometry, without reference to more general techniques in Riemannian geometry.

\subsection{Polyhedral approximation and the $n = 2$ case}

We will construct our examples using ``polyhedral approximation.'' For each integer $n \geq 2$ we will construct a convex $n$-polytope $\X$ (a closed subset of $\R^n$ with nonempty interior that is the intersection of finitely many closed half-spaces). Next we will glue two copies of $\X$ together via the identity map between their boundaries to obtain a topological manifold $\dbl\X$ that is called the \emph{double} of $\X$, and which is homeomorphic to $\Sp[n]$. Then we will construct a stable geodesic bouquet $G$ in a subspace of $\dbl\X$ that is isometric to a flat Riemannian manifold. $G$ will be composed of straight line segments. Finally, we will ``smooth'' $\dbl\X$ into a convex hypersurface of $\R^{n+1}$ with strictly positive curvature, while preserving the existence of a stable geodesic bouquet.

To illustrate the process by which we construct our examples, let us present the $n = 2$ case of our main result.

\begin{figure}[h]
    \centering
    \includegraphics[width=0.6\textwidth]{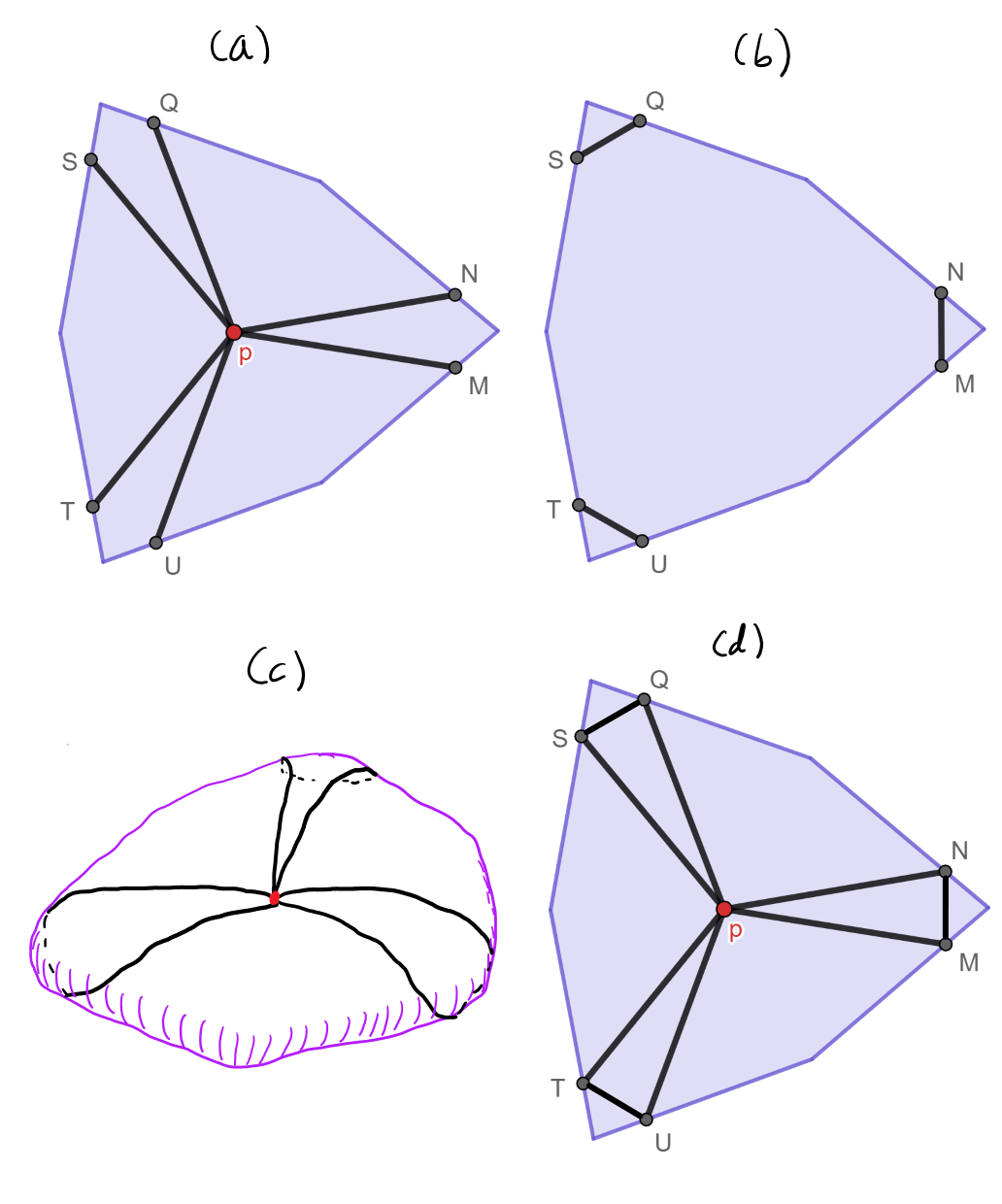}
    \caption{A ``trilateral bouquet'' in the double of a certain irregular hexagon with trilateral symmetry. The dark solid lines indicate segments of the geodesic bouquet. (a) Part of the bouquet in one copy of the hexagon in the double. (b) Part of the bouquet in the other copy of the hexagon in the double. (c) Stable geodesic bouquet in a ``smoothing'' of the double that has positive sectional curvature. (d) Billiard trajectories arise from combining (a) and (b) into a single figure.}
    \label{fig:TrilateralBouquet}
\end{figure}

\begin{family}[Trilateral Bouquet]
    \label{fam:TrilateralBouquet}
    Let $\X$ be the convex irregular hexagon with the symmetry group of an equilateral triangle and that is depicted in \cref{fig:TrilateralBouquet}(a). Up to dilation, it is determined by the interior angle $\theta$ closest to points $M$ and $N$, and we require that $\theta \in (\pi/3, \pi/2)$. Then $\dbl\X$ contains a stable geodesic bouquet $G$ with 3 loops, illustrated in \cref{fig:TrilateralBouquet}(a)--(b). $G$ is based at the center of the hexagon, and the rest is constructed so that the loops $pMN$, $pUT$ and $pSQ$ are \emph{billiard trajectories} in $\X$: that is, the line segments $pS$ and $SQ$ touch $\partial\X$ at the same angle at $S$, and the line segments $pQ$ and $SQ$ touch $\partial\X$ at the same angle at $Q$. The symmetry guarantees that $G$ is stationary. When $\theta < \pi/2$, $G$ avoids the vertices of the hexagon; when $\theta > \pi/3$, $G$ contains no closed geodesic. Later on we will prove that $G$ is stable in \cref{lem:TrilateralBouquetStable}.
    
    We can then smooth $\dbl\X$ into a hypersurface of $\R^3$ with positive curvature, as shown in \cref{fig:TrilateralBouquet}(c), in which there is a stable geodesic bouquet similar to $G$.\footnote{The proof of this assertion will be deferred the end of the paper, when we prove our main result, \cref{thm:StableNLoopsPosCurv}.}
\end{family}

\subsection{Billiards and the $n = 3$ case}

The information in \cref{fig:TrilateralBouquet}(a)--(b) can be summarized into the single diagram in \cref{fig:TrilateralBouquet}(d) by drawing all of the geodesic segments in the same convex polygon $\X$. This transforms the geodesics in the geodesic bouquet into \emph{billiard trajectories}. These are paths that begin at $p$, travel in straight lines reflect elastically off $\partial\X$, continue in straight lines, and so on until they eventually return to $p$. In this way, geodesics in $\dbl\X$ correspond to billiard trajectories in $\X$, with the understanding that every collision of the billiard trajectory corresponds to the geodesic passing from one copy of $\X$ in $\dbl\X$ to the other copy.

\begin{figure}[h]
    \centering
    \includegraphics[width=0.6\textwidth]{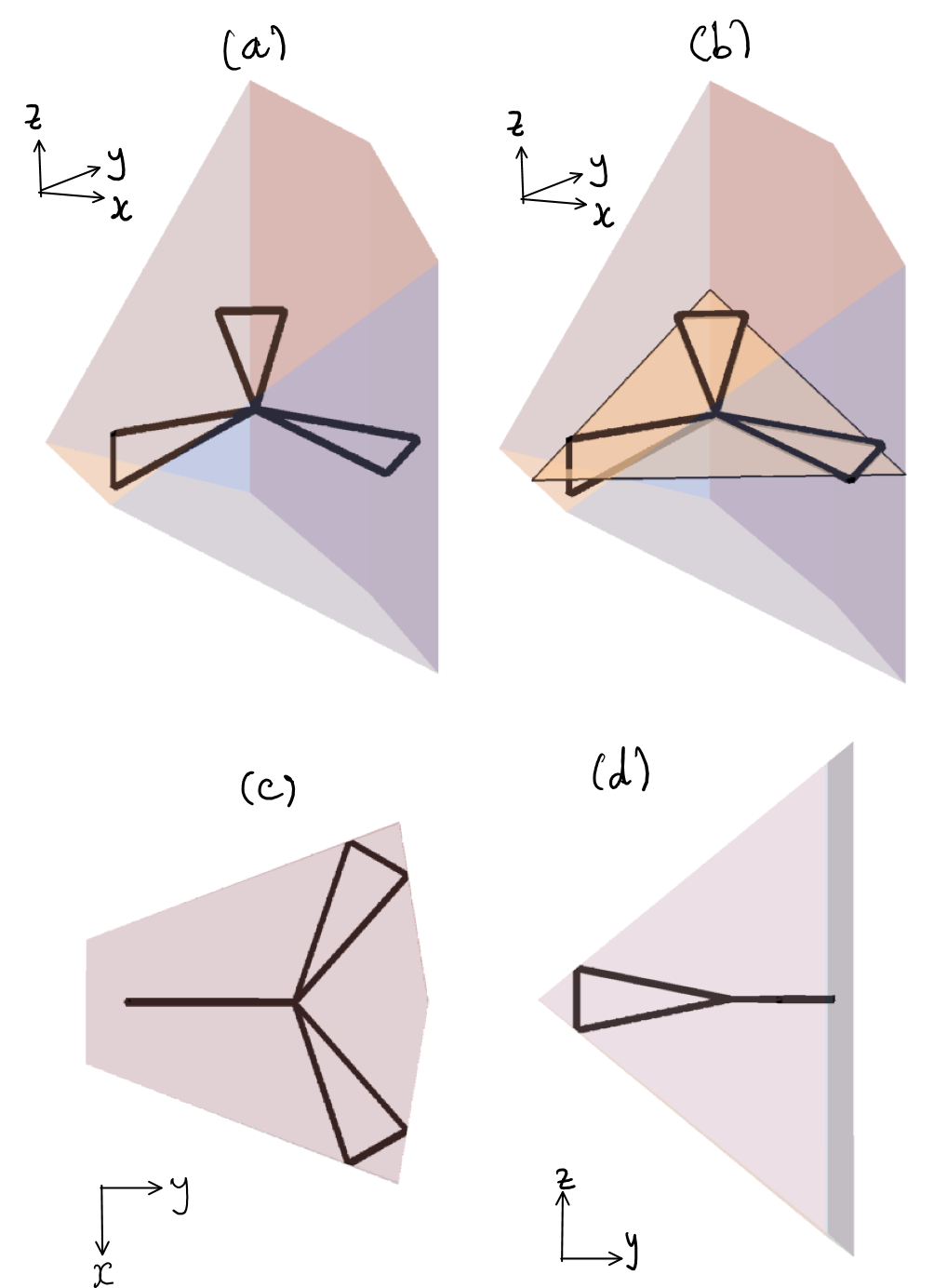}
    \caption{(a) The thick black lines depict three billiard trajectories (each forming a triangle) in a polyhedron $\X$ that correspond to a stable geodesic bouquet in $\dbl\X$. The trajectories all start and end at the same point. (b) Two of the billiard trajectories lie in the same plane as the triangle, while one does not. (c) and (d) show different views of $\X$ and the billiard trajectories.}
    \label{fig:BouquetFromTriangle}
\end{figure}

Accordingly, to prove the main result in a given dimension $n$, we will construct a convex $n$-dimensional polytope $\X$ and find billiard trajectories in $\X$ that will correspond to a stationary geodesic bouquet $G$ in $\dbl\X$. We will then prove that $G$ is stable, and prove that this stability is preserved after smoothing $\dbl\X$ to a convex hypersurface. The polytope and billiard trajectories for the $n = 3$ case are depicted in \cref{fig:BouquetFromTriangle}. The precise specifications of this construction, and the proof that it corresponds to a stable geodesic bouquet, will be presented in our proof of \cref{thm:StableNLoopsFlat}.

Observe from \cref{fig:BouquetFromTriangle}(b) that the three billiard trajectories do not lie in the same plane. This ``twisting'' of billiard trajectories relative to each other will lead to the stability of the corresponding geodesic bouquet, which will be proven in \cref{thm:StableNLoopsFlat}. Our constructions in higher dimensions will also hinge on finding configurations of billiard trajectories that are twisted relative to each other.

\subsection{Relation to the problem of capturing convex bodies}

After the completion of this work, Nabutovsky pointed out a connection between our results and the problem of ``capturing convex bodies with knots and links'', which has been discussed at some length on MathOverflow \cite{MO_CaptureSphereKnot,MO_CaptureConvexBodiesKnot}. Intuitively, a stable geodesic net in a convex hypersurface can in some cases be thought of as a net woven from string that cannot be stretched, and that ``captures'' the convex body that is bounded by that hypersurface. The meaning of capture is that the net cannot be slid off the convex body without stretching the strings. We caution that this analogy may not be exact if each edge of a stationary geodesic net $G$ is simply replaced by a single string and each vertex is formed just by gluing corresponding ends of string together. The reason is that the resulting ``string web'' would be physically taut when every sufficiently small perturbation must stretch at least one string, but that differs from the definition of a stable geodesic net: any sufficiently small perturbation of a stable geodesic net must increase the \emph{sum of edge lengths}.

To apply the analogy appropriately, one could model a stationary geodesic net as a mathematical link near each vertex in order to allow string to ``slide through'' the vertex. This would take into account variations that stretch some edges of a stationary geodesic net but contract others. A.~Geraschenko illustrated how to model a vertex of degree 3 as a link in \cite{MO_CaptureSphereKnot}. Under this implementation of the analogy, our \cref{fam:TrilateralBouquet} can be seen in retrospect as a modification of the net that captures an equilateral triangle, which was presented by A.~Petrunin in \cite{MO_CaptureConvexBodiesKnot}. Our modification produces a stable geodesic net that contains no closed geodesic.

\subsection{Organizational structure}

In \cref{sec:Definitions} we will formally define stable geodesic nets and other useful notions. In \cref{sec:PropertiesStableBouquets} we will derive some properties of stable geodesic bouquets, including some lower bounds on the index forms of geodesic loops. In \cref{sec:ConstructingBouquetInDouble} we will construct stable geodesic bouquets in the doubles of convex $n$-polytopes for each $n \geq 2$. We will use these constructions to prove our main result, \cref{thm:StableNLoopsPosCurv}.

\section{Definitions for Geodesic Bouquets}
\label{sec:Definitions}

Let $M$ be a Riemannian manifold. We will construct special stationary geodesic nets in $M$ called \emph{geodesic bouquets}, that are immersions of graphs $\bq_k$ which consist of a single vertex and $k$ loops at that vertex. Denote that vertex by $*$. Hence a stationary geodesic bouquet $G : \bq_k \to M$ has a \emph{basepoint} $p = G(*)$. $G$ is also composed of $k$ \emph{geodesic loops}, or geodesics that start and end at $p$.

For each integer $k \geq 1$, let $\imm{k}M$ be the set of piecewise $C^\infty$ immersions $\bq_k \to M$. Adapting \cite{Milnor_MorseTheory}, we endow this set with a topology induced by the metric where the distance between $G_1, G_2 \in \imm{k}M$ is defined as
\begin{equation}
    \max_{s \in \bq_k} \dist{G_1(s), G_2(s)} + \left(\int_{\bq_k} (\norm{G_1'(s)} - \norm{G_2'(s)})^2\,ds\right)^{1/2},
\end{equation}
where $G_i'(s)$ is the velocity vector of one of the loops of $G_i$, at a point $G(s)$, which is defined almost everywhere. The second summand is added to ensure that the length functional is a continuous function on $\imm{k}M$.

We then follow \cite{Staffa} and define an equivalence relation $\sim$ on the set of immersions $\bq_k \to M$ that relates immersions that are reparametrizations of each other. That is, $G_1 \sim G_2$ if $G_1 = \varphi \circ G_2$ for some homeomorphism $\varphi : \bq_k \to \bq_k$ that restricts to a diffeomorphism on each edge. We write the set of equivalence classes as $\qimm{k}M$, and endow it with the quotient topology. Denote the equivalence class of $G$ by $[G]$.

\subsection{Variations and stability of geodesic bouquets}

A \emph{variation} of a stationary geodesic bouquet $G$ is a homotopy $H : (-\epsilon, \epsilon) \times \Gamma \to M$ such that for each edge $e$ of $\Gamma$, $H|_{(-\epsilon, \epsilon) \times e}$ is a variation of the geodesic $H(0,e)$ in the usual sense. Sometimes we will write the variation as the one-parameter family of immersions $G_t = H(t,-)$. A \emph{vector field along $G$} is a continuous map $V : \Gamma \to TM$ that is piecewise $C^\infty$ on each edge of $\Gamma$ and satisfies $V(s) \in T_{G(s)}M$. We say that $V$ is \emph{tangent to $G$} if in addition we have $V(x) = 0$ for each vertex $x$ of $\Gamma$ with degree above 2, and for all other points $s$ that lie in the interior of an edge, $V(s)$ is either zero or tangent to the image of that edge.

For each variation $H$ of $G$, we write $\frac{dH}{dt}|_{t = 0}$ to mean the vector field along $G$ whose value at $s \in \Gamma$ is the velocity of the curve $t \mapsto H(t,s)$ at time $t = 0$. We also define the function $\length_H(t)$ whose value is the length of the image of the immersion $H(t,-)$. If $\frac{dH}{dt}|_{t = 0} = V$ then we say that $H$ \emph{is in the direction of $V$.}

We now define $G$ to be \emph{stable} if for every variation $H$ of $G$, either $\length_H''(0) > 0$ or $\frac{dH}{dt}|_{t = 0}$ is tangent to $G$. In this situation we call $G$ a \emph{stable geodesic net}.

It can be verified that if $G$ is stable, then $[G]$ is a strict local minimum in $\qimm{k}M$ under the length functional.

\subsection{Geodesic nets in polyhedral manifolds}

Consider the situation where $M$ is a \emph{polyhedral manifold} obtained by gluing together polyhedra. Formally, a polyhedral manifold is an $n$-dimensional triangulated topological manifold that is also a complete metric space, such that each simplex is isometric to an affine simplex in some Euclidean space \cite{Lebedeva_Smoothing3DPoly}. The complement of the $(n-2)$-skeleton of the triangulation is a (non-compact) flat Riemannian manifold which we call the \emph{smooth part} of $M$. All of the preceding definitions for geodesic nets would still make sense when $M$ is a polyhedral manifold, as long as the geodesic net lies in its smooth part. In particular, all edges of the geodesic net would locally be straight lines. Thus whenever we state that a polyhedral manifold $M$ contains a geodesic or geodesic net, we mean that the image of the geodesic or geodesic net is disjoint from the $(n-2)$-skeleton.

\section{Properties of Stable Geodesic Bouquets}
\label{sec:PropertiesStableBouquets}

Let $G$ be a stationary geodesic net that lies in the smooth part of a polyhedral manifold $M$. The goal of this section is to derive necessary conditions and sufficient conditions for the stability of $G$ that are easier to check than the definition of stability.

\subsection{Bouquet variations are constrained by their action at the basepoint}

It turns out that in a flat manifold, such as the smooth part of the double of a convex polytope, the variations $H$ of a stationary geodesic bouquet $G$ such that $\length_H''(0) \leq 0$ are severely constrained. In fact, those variations are in some sense completely determined by their behaviour at the basepoint of $G$. Hence we can reduce a lot of the analysis of such variations to linear algebra on the tangent space of the basepoint.

Let $\gamma : [0,1] \to M$ be a constant-speed geodesic loop of length $\ell$ and $V$ be a continuous and piecewise smooth vector field along $\gamma$. Let $H : (-\epsilon, \epsilon)_t \times [0,1]_s \to M$ be a variation of $\gamma$ in the direction of $V$. Let $V^\perp(t)$ denote the component of $V(t)$ orthogonal to $\gamma'(t)$. Then we have the second variation formula \cite[Theorem~6.1.1]{Jost_RiemannianGeometry}:
\begin{equation}
    \label{eq:2ndVariationGeodesicLoop}
    \length_H''(0) = \left.\ip{\frac{D}{dt}\frac{dH}{dt},\frac1\ell \gamma'}\right|_{(t,s) = (0,0)}^{(t,s) = (0,1)} + \frac1\ell \underbrace{\int_0^1\norm{\nabla_{\gamma'}V^\perp}^2 - \ip{R(\gamma',V^\perp)V^\perp,\gamma'}\,dt}_{Q_\gamma(V^\perp)},
\end{equation}
where $R$ is the Riemann curvature tensor of $M$, and where $Q_\gamma(-)$ is a quadratic form defined by $Q_\gamma(V) = \int_0^1\norm{\nabla_{\gamma'}V}^2 - \ip{R(\gamma',V)V,\gamma'}\,dt$.\footnote{Note that this definition differs from the usual definition of the index form of a unit-speed geodesic by a factor of $1/\ell$. However, this convention will simplify our subsequent formulas.}

For a stationary geodesic bouquet $G$ whose loops are the geodesics $\gamma_1, \dotsc, \gamma_k$, define its \emph{index form} to be the quadratic form
\begin{equation}
    \label{eq:IndexFormDef_GeodesicBouquet}
    Q_G(V) = \sum_{i = 1}^k Q_{\gamma_i}(V|_{\gamma_i}^\perp).
\end{equation}
Denote its nullspace by $\nullvar{G}$, that is, the space of vector fields $V$ along $G$ such that $Q_G(V) = 0$. By definition, $\nullvar{G}$ contains the vector fields tangent to $G$. Note that when the ambient manifold is flat, each summand of $Q_G$ becomes non-negative definite, so the same is true of $Q_G$.

\begin{lemma}
    \label{lem:IndexFormSecondVariation}
    Let $G$ be a stationary geodesic bouquet in a flat Riemannian manifold $M$. Let $V$ be a vector field along $G$. Then any variation $H$ of $G$ in the direction of $V$ satisfies $\length_H''(0) \geq 0$, where
    \begin{equation}
        \label{eq:2ndVariationGeodesicBouquet}
        \length_H''(0) = 0 \iff Q_G(V) = 0.
    \end{equation}
    Moreover, $G$ is stable if and only if $\nullvar{G}$ contains only vector fields tangent to $G$.
\end{lemma}
\begin{proof}
    Let $H_i : (-\varepsilon, \varepsilon) \times [0,1] \to M$ is the restriction of $H$ to $\gamma_i$. Consider the path $\alpha(t) = H(t,*) = H_i(t,0) = H_i(t,1)$. Thus if we let $\ell_i = \length(\gamma_i)$, then by \cref{eq:2ndVariationGeodesicLoop},
    \begin{align}
        \length_H''(0) 
        &= \sum_{i = 1}^k \left.\ip{\frac{D}{dt}\frac{dH_i}{dt}, \frac1{\ell_i} \gamma_i'}\right|_{(t,s) = (0,0)}^{(t,s) = (0,1)} + \frac1{\ell_i} \sum_{i = 1}^k Q_{\gamma_i}(V|_{\gamma_i}^\perp) \\
        &=  \underbrace{\ip{\left.\frac{D\alpha'}{dt}\right|_{t = 0}, \frac1{\ell_i} \sum_{i = 1}^k (\gamma_i'(1) - \gamma_i'(0))}}_0 + \frac1{\ell_i} \sum_{i = 1}^k Q_{\gamma_i}(V|_{\gamma_i}^\perp) \\
        &\geq 0,
    \end{align}
    where the first summand vanishes because of the stationarity condition in the definition of a geodesic bouquet. Evidently, $\length_H''(0)$ vanishes if and only if each $Q_{\gamma_i}(V|_{\gamma_i}^\perp)$ also vanishes, which is equivalent to the vanishing of $Q_G(V)$.
    
    Suppose that $G$ is stable. Then for any $V \in \nullvar{G}$, the definition of stability forces $V$ to be tangent to $G$. Conversely, suppose that every vector field in $\nullvar{G}$ is tangent to $G$. Then for any variation $H(t,s)$ of $G$ such that $\length_H''(0) = 0$, $Q_G(V) = 0$ where $V = \frac{dH}{dt}|_{t = 0}$. So $V \in \nullvar{G}$ and must be tangent to $G$. Hence $G$ is stable.
\end{proof}

In order to prove the stability of the stationary geodesic nets that we will construct later, we will need some lower bounds on the index forms. In particular we will prove lower bounds in terms of quantities and conditions that can be computed or checked at the basepoint.

The following proposition proves a lower bound that applies even for manifolds with small but positive sectional curvature. This helps us to apply it later to parts of the smoothings of doubles of convex polytopes that have small curvature.

Let $G$ be a stationary geodesic bouquet in $M$, one of whose loops is $\gamma : [0,1] \to M$. If $M$ is flat, then $Q_G(W) = 0$ only if $Q_\gamma(W|_\gamma^\perp) = 0$, which can only happen if $W|_\gamma^\perp$ is parallel along $\gamma$. For this reason, it is natural to ask whether each $v \in T_pM$ can be extended to a vector field $V$ along $\gamma$ such that $V^\perp$ is parallel along $\gamma$. Such an extension will be impossible for some vectors $v \in T_pM$, and we can measure the obstruction to this extension using the \emph{parallel defect operator} $\opd{\gamma} : T_pM \to T_pM$ defined as follows.

\begin{definition}[Parallel defect operator, kernel]
    \label{def:ParallelDefectOperator}
    Let $\gamma : [0,1] \to M$ be a geodesic loop with parallel transport map $P : T_pM \to T_pM$. For each $t \in [0,1]$, let $\pi_t : T_{\gamma(t)}M \to T_{\gamma(t)}M$ be the projection onto the orthogonal complement of $\gamma'(t)$. Then the \emph{parallel defect operator of $\gamma$} is $\opd\gamma = P\pi_0 - \pi_1$. The \emph{parallel defect kernel} of $\gamma$ is $\ker\opd\gamma$.
\end{definition}

The value $\norm{\opd\gamma v}$ quantifies how difficult it is to extend $v$ in the aforementioned manner. Roughly speaking, it turns out that if such an extension is difficult for certain $v$, then any extension of $v$ into a vector field $V$ along $\gamma$ gives a large value for $Q_\gamma(V^\perp)$. The following proposition makes this precise by bounding $Q_\gamma(V)$ from below in terms of $\norm{\opd\gamma v}$ in flat manifolds.

\begin{proposition}
    \label{prop:IndexFormLowerBound_ParallelDefect_Flat}
    Let $M$ be an $n$-dimensional flat Riemannian manifold. Suppose that $M$ contains a geodesic loop $\gamma : [0,1] \to M$ of length $\ell$. Then for any vector field $V$ along $\gamma$,
    \begin{equation}
        \label{eq:IndexFormLowerBound_GeodesicLoop_Flat}
        Q_\gamma(V^\perp) \geq \norm{\opd\gamma V(0)}^2.
    \end{equation}
\end{proposition}
\begin{proof}
    Let $F_1, \dotsc, F_{n-1},\gamma'$ be a parallel orthonormal frame along $\gamma$. Let $V^\perp = \sum_{i = 1}^{n-1}c_iF_i$ for piecewise smooth functions $c_i : [0,1] \to \R$. To bound the index form, observe that $\int_0^1 \norm{\nabla_{\gamma'} V^\perp}^2\,dt = \int_0^1 \sum_{i = 1}^{n-1}c_i'(t)^2\,dt$ is twice of the energy of the path $\alpha(t) = (c_1(t),\dotsc, c_{n-1}(t))$ in $\R^{n-1}$ \cite[eq.~(1.4.8)]{Jost_RiemannianGeometry}. If we fix the endpoints of this path, then the energy is minimized by the constant-speed straight path $\overline{\alpha}$ with the same endpoints. This minimum energy is $\frac12\sum_{i = 1}^{n-1} (c_i(b) - c_i(a))^2$. Define $P$, $\pi_0$ and $\pi_1$ as in \cref{def:ParallelDefectOperator}. Then
\begin{align*}
    \int_0^1 \norm{\nabla_{\gamma'}V^\perp}^2 \,dt
    &\geq \sum_{i = 1}^{n-1} (c_i(1) - c_i(0))^2 \\
    &= \norm{\sum_{i = 1}^{n-1} (c_i(0) - c_i(1))F_i(1)}^2 \\
    &= \norm{\sum_{i = 1}^{n-1} c_i(0)F_i(1) - \sum_{i = 1}^{n-1} c_i(1)F_i(1)}^2 \\
    &= \norm{P\sum_{i = 1}^{n-1} c_i(0)F_i(0) - \sum_{i = 1}^{n-1} c_i(1)F_i(1)}^2 \\
    &= \norm{P\pi_0V(0) - \pi_1 V(1)}^2 \\
    (V(0) = V(1)) &= \norm{\opd\gamma V(0)}^2. \qedhere
\end{align*}
\end{proof}

As a consequence, parallel defect kernels are related to the stability of geodesic bouquets, as expressed in the following result.

\begin{corollary}
    \label{cor:NullvarsToIntersectKernels_Flat}
    Let $G$ be a stationary geodesic bouquet in a flat manifold $M$, based at $p$, with loops $\gamma_1, \dotsc, \gamma_k$. Then $\setb{V(*)}{V \in \nullvar{G}} = \bigcap_{i = 1}^k \ker\opd{\gamma_i}$ and $G$ is stable if and only if $\bigcap_{i = 1}^k \ker\opd{\gamma_i} = \{0\}$.
\end{corollary}
\begin{proof}
    $M$ is flat, so $Q_G$ and each $Q_{\gamma_i}$ are non-negative definite. ($\subset$): Let $V \in \nullvar{G}$, so $Q_G(V) = 0$. By \cref{eq:IndexFormDef_GeodesicBouquet} and \cref{prop:IndexFormLowerBound_ParallelDefect_Flat}, this implies that $v \in \ker\opd{\gamma_i}$ for all $1 \leq i \leq k$, where $v$ is the value of $V$ at the basepoint of $G$.
    
    ($\supset$): let $v \in \bigcap_i\opv{\gamma_i}$. For each $i$, let $w_i$ be the component of $v$ that is orthogonal to $\gamma_i'(0)$. Extend $w_i$ to a parallel vector field $W_i$ along $\gamma_i$ by parallel transport. Now define
    \begin{equation}
        V_i(t) = W_i(t) + \left((1 - t)\ip{\gamma_i'(0),v} + t\ip{\gamma_i'(1),v}\right)\gamma_i'(t).
    \end{equation}
    Piecing all the $V_i$'s together gives a vector field in $\nullvar{G}$ that evaluates to $v$ at $p$.
    
    Now we prove the final statement in the corollary. By \cref{lem:IndexFormSecondVariation}, $G$ is stable if and only if $\nullvar{G}$ contains only vector fields tangent to $G$. If $\nullvar{G}$ contains only vector fields tangent to $G$, then in particular all of those vector fields must vanish at the basepoint, so as shown above, $\bigcap_i\opv{\gamma_i} = \{0\}$. Conversely, if $\bigcap_i\opv{\gamma_i} = \{0\}$, then every vector field $V \in \nullvar{G}$ must vanish at the basepoint. However, since $Q_G(V) = 0$, we also have $Q_\gamma(V|_\gamma^\perp) = 0$ for each loop $\gamma$ of $G$. This means that $V|_\gamma^\perp$ is parallel along $\gamma$. However, $V|_\gamma^\perp$ vanishes at the basepoint, so in fact it must be identically zero. That is, $V|_\gamma$ is tangent to $\gamma$. To sum up, $V$ must be tangent to $G$.
\end{proof}

\section{Constructing Stable Geodesic Bouquets}
\label{sec:ConstructingBouquetInDouble}

The goal of this section is to prove our main result, the bulk of which is a consequence of the following theorem:

\begin{restatable}{theorem}{ThmStableNLoopsFlat}
    \label{thm:StableNLoopsFlat}
    For every integer $n \geq 3$, there exists a compact convex $n$-polytope $\X$ such that $\dbl\X$ contains a stable geodesic bouquet with $n$ loops. For $n = 2$, there exists a hexagon $\X$ such that $\dbl\X$ contains a stable geodesic bouquet with 3 loops. Moreover, in each of these stable geodesic bouquets, no two tangent vectors of the loops at the basepoint are parallel.
\end{restatable}

This theorem would almost imply our main result, as long as one believes that the stability of geodesic bouquets is preserved under slight smoothings of $\dbl\X$. The $n = 2$ case of \cref{thm:StableNLoopsFlat} will be settled by \cref{fam:TrilateralBouquet} from the Introduction. As shown in \cref{fig:TrilateralBouquet}(d), the loops $pMN$, $pUT$ and $pSQ$ are billiard trajectories in the hexagon $\X$. There is indeed a correspondence between billiard trajectories in a convex polytope $\X$ and its double $\dbl\X$, which we will soon explain. To prove \cref{thm:StableNLoopsFlat} for $n \geq 3$, we will construct billiard trajectories in some convex $n$-polytope $\X$ which, under that correspondence, will correspond to the geodesic loops of the desired stable geodesic bouquet.

Hence we will start this section with some background and notation for billiards and formalize this correspondence.

\subsection{Definitions for polytopes and billiards}
\label{sec:DefnPolytopesBilliards}

Let $\X$ be a convex $n$-polytope. When $\X$ is compact, $\partial\X$ is homeomorphic to $\Sp[n-1]$ and has a natural cellular decomposition such that each $m$-cell (for $m \geq 1$) is isometric to an $m$-polytope. We require this decomposition to be ``irreducible'' in the following sense: no two distinct $m$-cells can lie in the same $m$-dimensional affine subspace of $\R^n$. The 0-cells are called the \emph{vertices} of $\X$, while the $(n-1)$-cells are called the \emph{faces}.

A \emph{supporting hyperplane} of $\X$ is a hyperplane $H \subset \R^n$ that intersects $\X$ at some face $F$. In this situation, $H$ is called the \emph{supporting hyperplane} of $F$, and the \emph{supporting half-space} of $F$ is a half-space $\HS$ of $\R^n$ containing $\X$ whose boundary is the supporting hyperplane of $F$. We say that $\HS$ is a \emph{supporting half-space} of $\X$.

The interior and boundary of a subspace $A$ of a topological space $S$ is denoted by $\itr{A}$ and $\partial A$ respectively.

Let $\X$ be a convex $n$-polytope. When we speak of the \emph{boundary} of an $m$-cell $A$ for $m \geq 1$, which by an abuse of notation we will denote by $\partial A$, we mean the union of the $(m-1)$-cells that are contained in $A$ (note that $A$ is homeomorphic to a closed disk). By the interior of $A$ we mean $A \setminus \partial A$, and denote it by $\itr{A}$.

If $F$ is a face of $\X$, let $R_F : \R^n \to \R^n$ denote the reflection about the supporting hyperplane of face $F$.

A \emph{billiard trajectory} in a convex $n$-polytope $\X$ is a sequence of points $x_0, x_1, \dotsc, x_k, x_{k+1} \in \X$ such that for each $1 \leq i \leq k$, $x_i$ lies in the interior of some face $F_i$ of $\X$ such that $F_i \neq F_{i+1}$ for $1 \leq i \leq k-1$. We also require that for $0 \leq i \leq k$, the unit vectors $u_i = (x_{i+1} - x_i)/\norm{x_{i+1} - x_i}$ must satisfy $u_{i+1} = dR_{F_i}u_i$, where $dR_{F_i}$ is the differential (linear part) of the affine transformation $R_{F_i}$. (That is, if $R_{F_i}x = Ax + b$ for a matrix $A$ and $b \in \R^n$, then $dR_{F_i} = A$.) The points $x_1, \dotsc, x_k$ are called the \emph{collisions} of the billiard trajectory, and we say that the billiard trajectory collides with faces $F_1, \dotsc, F_k$ in that order. Such a billiard trajectory is \emph{proper} if $x_0, x_{k+1} \in \itr\X$. We will usually represent billiard trajectories as a path $\beta : [a,b] \to \X$, parametrized at constant speed, consisting of line segments (called the \emph{segments} of $\beta$) joining $\beta(a) = x_0$ to $x_1$, $x_1$ to $x_2$ and so on until it joins $x_k$ to $x_{k+1} = \beta(b)$. Clearly, this is equivalent to the representation in terms of a sequence of points. $\beta$ is a \emph{billiard loop} if it is proper and $\beta(a) = \beta(b)$. We say that $\beta$ is \emph{periodic} if it is a billiard loop satisfying $\beta'(a) = \beta'(b)$ (see \cref{fig:GeodesicBilliardCorrespondence}(a) for an illustration of a periodic billiard trajectory).

\subsection{Correspondence between geodesics and billiards in the double of a convex polytope}

By the \emph{double} of a compact topological manifold with boundary $M$ we mean the closed topological manifold that is the quotient $\dbl{M} = M \cup_\partial M = (M \times \{0,1\})/\sim$ where $(x,0) \sim (x,1)$ for all $x \in \partial M$. It comes with a natural quotient map $\pi : \dbl M \to M$ that sends all $(x,i)$ to $x$.

Given a compact convex $n$-polytope $\X$, its double $\dbl\X$ is homeomorphic to $\Sp[n]$, and it inherits a natural cellular decomposition from $\X$. The $(n-2)$-skeleton of $\dbl\X$ is called its \emph{singular set}, denoted by $\dbl[sing]\X$. The complement of this singular set in $\dbl\X$ is a non-compact flat manifold, denoted by $\dbl[sm]\X$.

\begin{figure}[h]
    \centering
    \includegraphics[height=5cm]{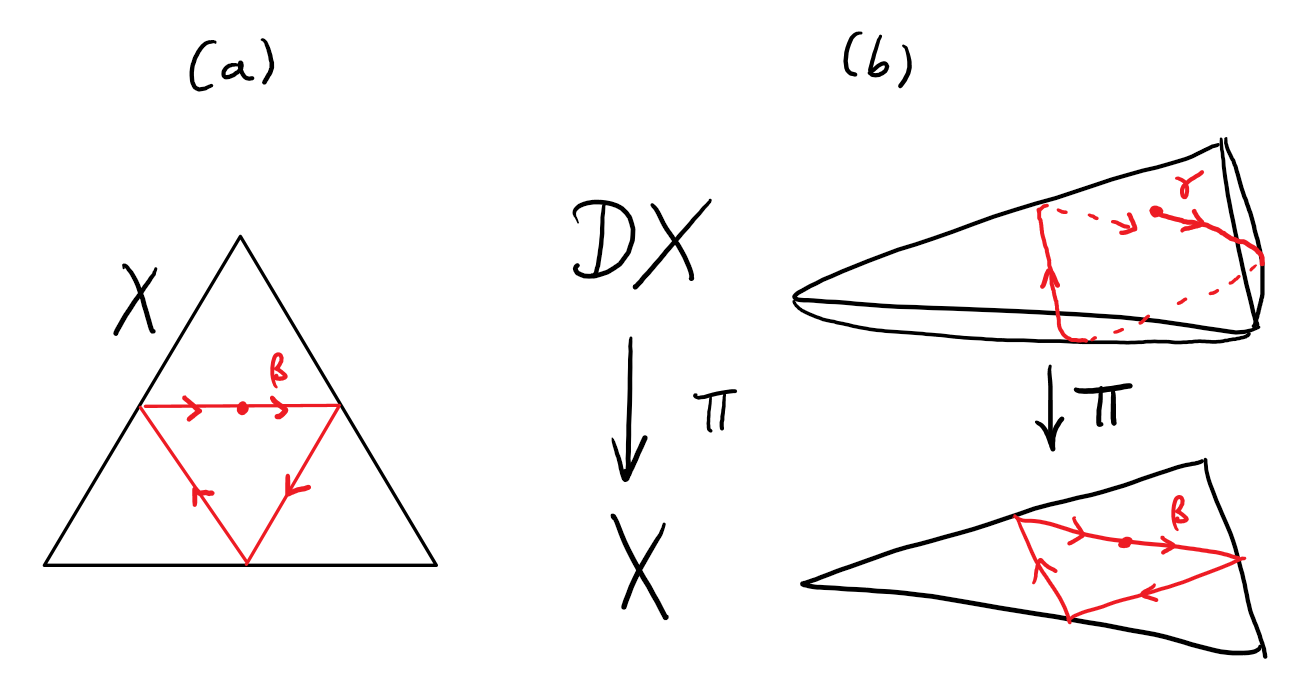}
    \caption{(a) A periodic billiard trajectory that makes 3 collisions inside an equilateral triangle $\X$. (b) The corresponding geodesic $\gamma$ in $\dbl\X$. Dashed lines indicate parts of $\gamma$ that lie on the ``bottom'' copy of $\X$ in $\dbl\X$. $\gamma$ starts on the ``top'' copy of $\X$ and ends on the ``bottom'' copy, so $\gamma$ is not closed.}
    \label{fig:GeodesicBilliardCorrespondence}
\end{figure}

Let $\gamma : [a,b] \to \dbl\X$ be a geodesic that avoids the singular set. (Henceforth we will just say that $\gamma$ is a geodesic in $\dbl\X$, and leave it implicitly understood that it avoids the singular set.) Then the composition $\beta : [a,b] \xrightarrow{\gamma} \dbl\X \xrightarrow{\pi} \X$ is a billiard trajectory (see \cref{fig:GeodesicBilliardCorrespondence}(b)). Conversely, given a billiard trajectory $\beta : [a,b] \to \X$, and a choice of $x \in \pi^{-1}(\beta(a))$, there is a unique geodesic $\gamma : [a,b] \to \dbl$ such that $\gamma(a) = x$ and $\pi \circ \gamma = \beta$. Hence there is this correspondence between geodesics in $\dbl\X$ and billiard trajectories in $\X$. Under this correspondence, closed geodesics in $\dbl\X$ that start and end in the interior of an $n$-cell $A$ correspond to periodic billiard trajectories that have an even number of collisions, because closed geodesics have to exit and enter $A$ an even number of times in total. If a periodic billiard trajectory has an odd number of collisions in $\X$, then its corresponding geodesic in $\dbl\X$ is not closed (see \cref{fig:GeodesicBilliardCorrespondence}(b)).

\subsection{Geodesic bouquets from billiard loops that collide twice}

If we are trying to prove \cref{thm:StableNLoopsFlat} for $n \geq 3$, then our geodesic bouquet is allowed to have $n$ loops, where the bouquet is inside $\dbl\X$ for some convex $n$-polytope $\X$. The geodesic loops $\gamma_i$ in the bouquet will correspond to billiard trajectories $\beta_i$ in $\X$ that have a particularly simple form: they only collide twice with $\partial\X$. This allows $\opv{\gamma_i}$ to be computed explicitly as follows.

\begin{lemma}
    \label{lem:2CollisionBilliadLoopOPV}
    Let $\X$ be a convex $n$-polytope and $\gamma : [0,1] \to \dbl\X$ be a geodesic loop based at $p$ such that $\beta = \pi \circ \gamma$ is a proper billiard trajectory that makes two collisions. Then $\opv\gamma$ is the hyperplane in $T_p\dbl\X$ that is orthogonal to $\gamma'(0) + \gamma'(1)$.
\end{lemma}
\begin{proof}
    Note that the three segments of $\beta$ form a triangle and therefore lie in the same plane $W_\beta$. The definition of a billiard trajectory implies, by some elementary geometry, that $W_\beta$ must contain the normal vectors of the faces that $\beta$ collides with. Consequently, $\X$ locally looks like the product $(W_\beta \cap \X) \times \R^{n-2}$ near $\beta$. This in turn implies that $\dbl\X$ locally looks like $\dbl(W_\beta \cap \X) \times \R^{n-2}$ near $\gamma$. This constrains the parallel transport map $P$ of $\gamma$: if we write $W_\gamma = \vspan\{\gamma'(0), \gamma'(1)\}$, then $P$ must fix $W_\gamma^\perp$. Since $\dbl\X$ is orientable, $\det P = 1$ so $P$ must act on $W_\gamma$ as the rotation that brings $\gamma'(0)$ to $\gamma'(1)$.
    
    Consider the action of $\opd\gamma$ on a basis consisting of $\gamma'(0) + \gamma'(1)$, $\gamma'(0) - \gamma'(1)$, and vectors from $W_\gamma^\perp$. Note that $W_\gamma^\perp$ is an invariant subspace of the linear operators $P$, $\pi_0$ and $\pi_1$ in the definition of $\opd\gamma$. This implies that $W_\gamma^\perp \subset \ker\opd\gamma$. As for the vector $\gamma'(0) - \gamma'(1)$, it projects via $\pi$ to the internal angle bisector (within the plane $W_\beta$) of the angle of the triangle $\beta$ at $\pi(p)$. Thus $v_0 = \pi_0(\gamma'(0) - \gamma'(1))$ and $v_1 = \pi_1(\gamma'(0) - \gamma'(1))$ are vectors of the same length. It can also be verified that $(v_0,\gamma'(0))$ and $(v_1,\gamma'(1))$ are orthogonal bases of $W_\beta$ with the same orientation, so $Pv_0 = v_1$ and $D_\gamma(\gamma'(0) - \gamma'(1)) = 0$.

    Therefore $\ker\opd\gamma$ contains the hyperplane $W_\gamma^\perp \oplus \vspan\{\gamma'(0) - \gamma'(1)\}$, and it remains to prove that $\opd\gamma(\gamma'(0) + \gamma'(1)) \neq 0$. To prove this, observe that $\gamma'(0) + \gamma'(1)$ projects via $\pi$ to the external angle bisector of the angle of $\beta$ (within $W_\beta$) at $\pi(p)$. This implies that $\pi_0(\gamma'(0) + \gamma'(1)) = -\lambda v_0$ and $\pi_1(\gamma'(0) + \gamma'(1)) = \lambda v_1$ for some $\lambda > 0$. Hence $\opd\gamma(\gamma'(0) + \gamma'(1)) = P(-\lambda v_0) - \lambda v_1 = -2\lambda v_1 \neq 0$.
\end{proof}

\begin{lemma}
    \label{lem:TrilateralBouquetStable}
    The geodesic bouquets in \cref{fam:TrilateralBouquet} are stable.
\end{lemma}
\begin{proof}
    By \cref{lem:2CollisionBilliadLoopOPV}, the $\opv\gamma$'s bisect the angles $\angle MpN$, $\angle UpT$ and $\angle SpQ$. Their intersection has dimension 0. Apply \cref{cor:NullvarsToIntersectKernels_Flat}.
\end{proof}

\begin{remark}
    Every stable geodesic bouquet in the double of a convex $n$-polytope derived only from billiard trajectories that collide twice must have at least $n$ loops, by the dimension formula. Conversely, $n$ loops is enough when $n \geq 3$.
\end{remark}

Before we construct $\X$ and the stationary geodesic bouquet in $\dbl\X$, we need the following linear-algebraic lemma, to be proven in Appendix~\ref{sec:LinearAlgebraLemma}.

\begin{restatable}{lemma}{LemmaGenPosHyperplanesFromSimplex}
    \label{lem:GenPosHyperplanesFromSimplex}
    Let $X$ be an $n$-dimensional vector space. Let $x_1, \dotsc, x_n \in X$ be nonzero vectors that span an $(n-1)$-dimensional subspace, satisfying the property that any $n-1$ of them are linearly independent. Then there exist $(n-1)$-dimensional subspaces of $X$, $\Pi_1,\dotsc, \Pi_n$ such that $x_i \in \Pi_i$ for all $i$ and $\bigcap_{i=1}^n \Pi_i = \{0\}$.
\end{restatable}

Now we are ready to prove \cref{thm:StableNLoopsFlat}, restated as follows.

\ThmStableNLoopsFlat*

\begin{proof}
    The dimension 2 case has been handled by \cref{lem:TrilateralBouquetStable}. Let $n \geq 3$ be an integer. Let $\sigma$ be a regular $(n-1)$-simplex in $\R^{n-1}$ with vertices $\bar{x}_1, \dotsc, \bar{x}_n$ on the unit sphere, located so that its barycenter is at the origin $o \in \R^{n-1}$. It satisfies the following key properties, where $o \in \R^{n-1}$ denotes the origin.
    \begin{enumerate}
        \item $\sigma$ is non-degenerate. That is, its vertices are affinely independent, i.e. from any vertex, the $n-1$ vectors to the other vertices are linearly independent.
        
        \item $\bar{x}_1 + \dotsb + \bar{x}_n = 0$.
        
        \item For each $1 \leq i < j \leq n$, the angle $\angle \bar{x}_i\bar{x}_jo < \pi/4$.\footnote{This is a consequence of basic Euclidean geometry.}
    \end{enumerate}
    
    $\sigma$ has a ``geometric dual'', another regular simplex $\hat\sigma$ whose faces are tangent to the unit sphere at the points $\bar{x}_i$ (see \cref{fig:GeodesicBouquetFromSimplex_SetupNotation}(a)).
    
    Now consider the unit vectors $x_i = (\bar{x}_i,0) \in \R^n$, which also must be affinely independent. The convex $n$-polytope $\X_0 := \hat\sigma \times \R$ is non-compact (see \cref{fig:GeodesicBouquetFromSimplex_SetupNotation}(b)), but we will modify it along a 1-parameter family of convex polytopes $\{\X_\theta\}_{\theta \in [0,\pi/2)}$ which for all $\theta \in (0, \pi/4]$ are compact and contain all of the points $x_i$. We will then pick some value of $\theta$ slightly larger than $\pi/4$ and find billiard trajectories in $\X_\theta$ that are close to the line segments $ox_i$, and show that these billiard trajectories correspond to the desired geodesic bouquet in $\dbl\X_\theta$. More specifically, each $x_i$ lies on a unique face of $\X_0$, which has supporting half-space $\HS_i^0$ (see \cref{fig:GeodesicBouquetFromSimplex_SetupNotation}(b)). We will find a 1-parameter family of half-spaces $\HS_i^\theta$ whose boundaries contain $x_i$ and are ``rotations of $\HS_i^0$ by angle $\theta$''. Then we will replace each $\HS_i^0$ with the convex set $\HS_i^{-\theta} \cap \HS_i^\theta$ to get the convex $n$-polytope $\X_\theta = \bigcap_{i = 1}^n (\HS_i^{-\theta} \cap \HS_i^\theta)$.
    
    The directions of ``rotation'' will be determined by vector subspaces of $\R^n$ that are obtained from \cref{lem:GenPosHyperplanesFromSimplex}, and which will eventually correspond to the kernels of some parallel defect operators. To apply the lemma we must check that any $n-1$ of the vectors $x_1, \dotsc, x_n$ have to be linear independent. Suppose the contrary, that without loss of generality, $a_1x_1 + \dotsb + a_{n-1}x_{n-1} = 0$ for some $a_i \in \R$. We will deduce a contradiction with the affine independence of $x_1,\dotsc,x_n$. For any $b_i \in \R$ we have
    \begin{equation}
        \sum_{i = 1}^{n-1} b_i(x_i - x_n) 
        = \sum_{i = 1}^{n-1} b_ix_i - \left(\sum_{i = 1}^{n-1} b_i\right) x_n
        = \sum_{i = 1}^{n-1} b_ix_i + \left(\sum_{i = 1}^{n-1} b_i\right)\left(\sum_{i = 1}^{n-1} x_i\right) = \sum_{i = 1}^{n-1} c_ix_i,
    \end{equation}
    where $c_i = b_i + \sum_{i = 1}^{n-1} b_i$. If we can choose the $b_i$'s such that $c_i = a_i$ then we would have a contradiction against the affine independence of $x_1,\dotsc,x_n$. It is possible to choose it this way, because the vector $(c_1,\dotsc,c_{n-1})$ can be obtained by multiplying the vector $(b_1,\dotsc,b_{n-1})$ by a matrix whose diagonal entries are 2 and all other entries are 1. Gaussian elimination reveals that this matrix is invertible and thus we can solve for the $b_i$'s after setting $c_i = a_i$. 
    
    Therefore any $n-1$ of the vectors $x_1, \dotsc, x_n$ must be linearly independent. Since the vectors sum to 0, they must span a hyperplane. Using \cref{lem:GenPosHyperplanesFromSimplex}, choose hyperplanes $\Pi_i$ passing through the origin and $x_i$ such that $\bigcap_i \Pi_i = \{0\}$ (see \cref{fig:GeodesicBouquetFromSimplex_SetupNotation}(c)).
    
    \begin{figure}[h]
    \centering
    \includegraphics[width=\textwidth]{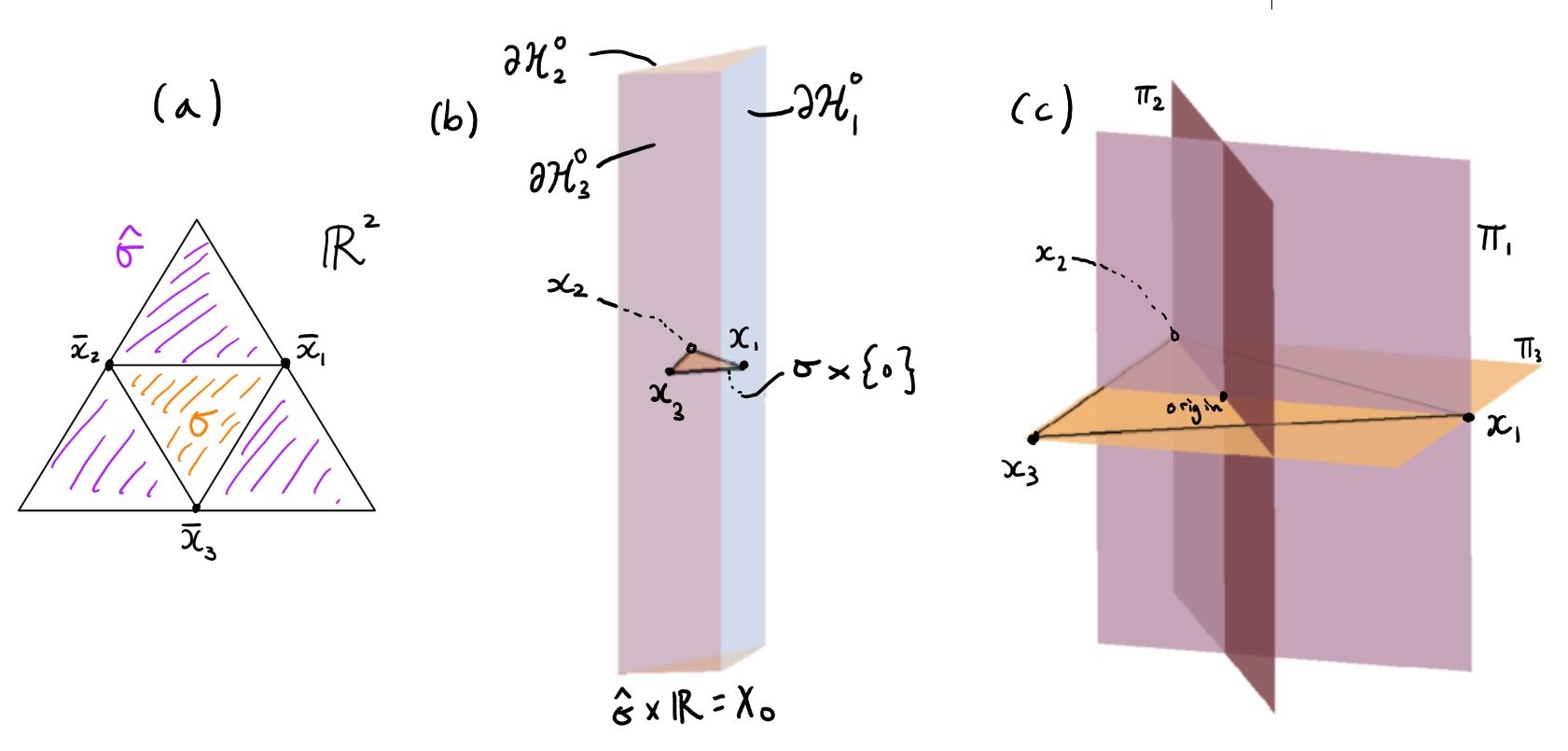}

    \caption{An illustration of key objects in the proof of \cref{thm:StableNLoopsFlat} for $n = 3$. (a) $\sigma$ is a regular 2-simplex in $\R^2$ that has a ``geometric dual'' $\hat\sigma$. (b) The non-compact 3-polytope that we will modify into the polytope we need. A face of the polytope is labeled $\partial \HS_i^0$ to indicate that its supporting half-space is $\HS_i^0$. (c) The $\Pi_i$'s are hyperplanes derived from \cref{lem:GenPosHyperplanesFromSimplex} that will eventually be the kernels of parallel defect operators.}
    \label{fig:GeodesicBouquetFromSimplex_SetupNotation}
\end{figure}
    
    Let $\theta \in (-\pi/2,\pi/2)$ be any angle. For each $1 \leq i \leq n$, let the plane $W_i$ be spanned by $x_i$ and the normal vector of $\Pi_i$. Let $\HS_i^\theta$ be the image of $\HS_i^0$ under the rotation that acts as the identity on $W_i^\perp + x_i$ ($W_i^\perp$ translated by vector $x_i$) and as the rotation by angle $\theta$ on $W_i$. This requires choosing an orientation on $W_i$, but we will define $\X_\theta = \bigcap_{i = 1}^n (\HS_i^{-\theta} \cap \HS_i^\theta)$, and $\X_\theta$ will not depend on this choice.

    Note that $\X_\theta \subset \X_0$ for all $\theta \geq 0$. Now we show that $\X_\theta$ is compact for $\theta > 0$. This is because the construction in \cref{lem:GenPosHyperplanesFromSimplex} actually implies that $\Pi_n = \R^{n-1} \times \{0\}$. The only way $\X_\theta$ could be non-compact is if it contains a ``vertical ray'', that is, a ray orthogonal to $\Pi_n$. However, such a ray would have to intersect $\partial(\HS_n^{-\theta} \cap \HS_n^\theta)$. Therefore $\X \cap (\HS_n^{-\theta} \cap \HS_n^\theta)$ would already be compact, thus $\X_\theta$ is also compact.
    
\begin{figure}[p]
    \centering
    \includegraphics[width=0.9\textwidth]{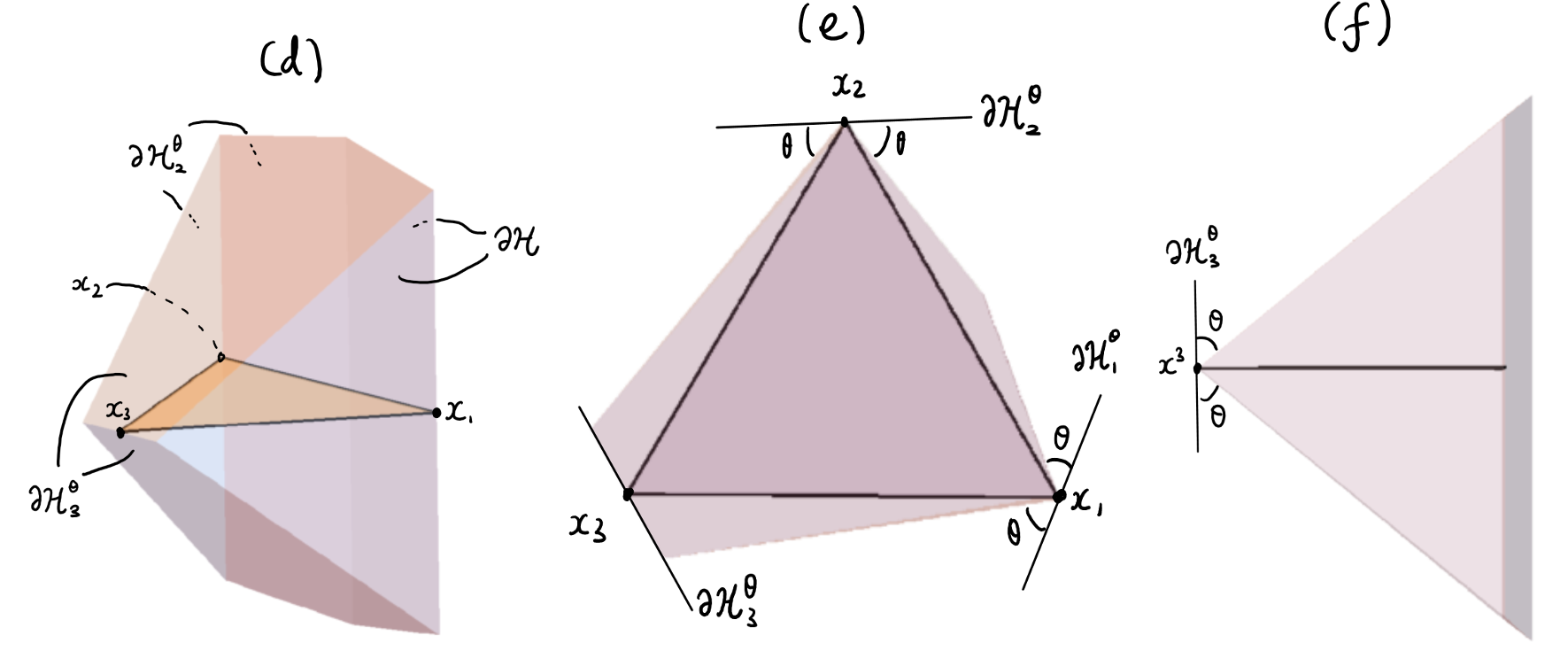}
    \includegraphics[width=0.9\textwidth]{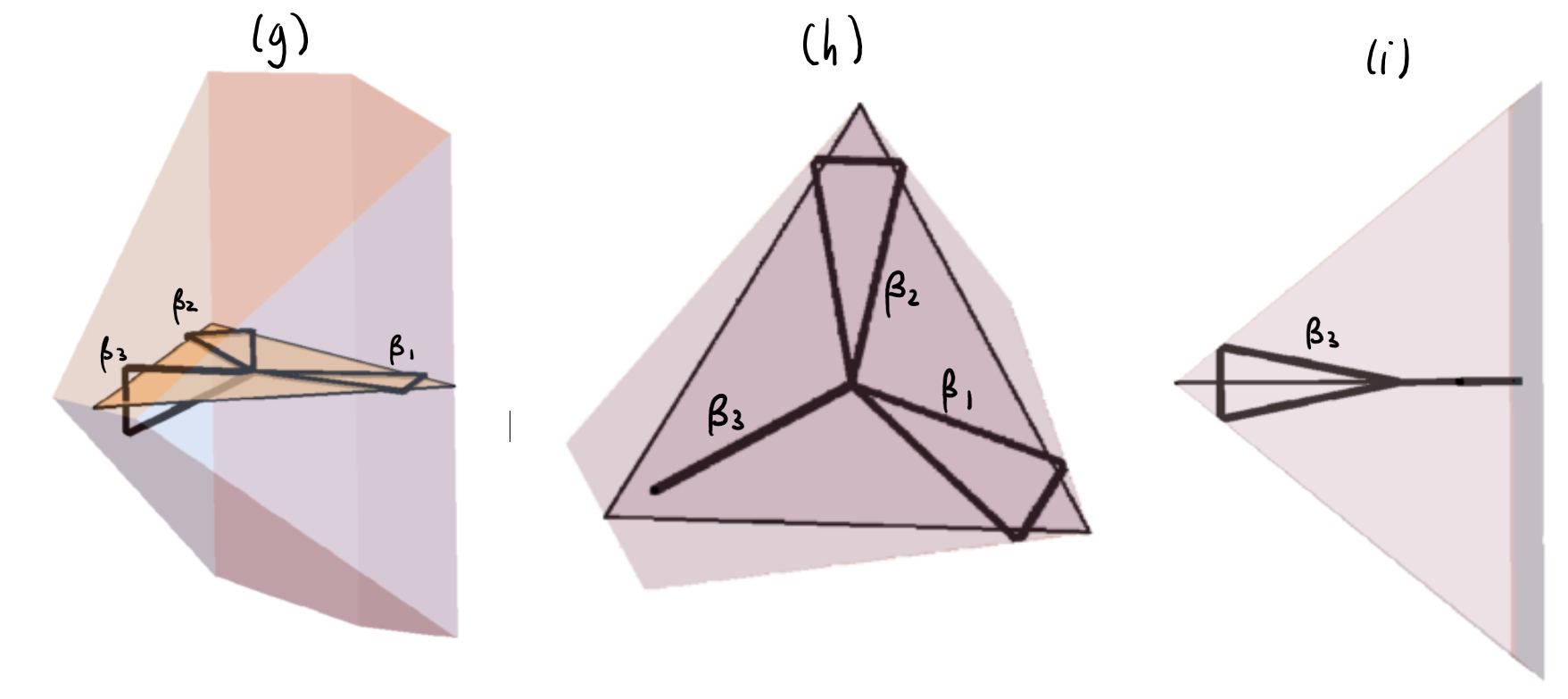}
    \includegraphics[scale=0.6]{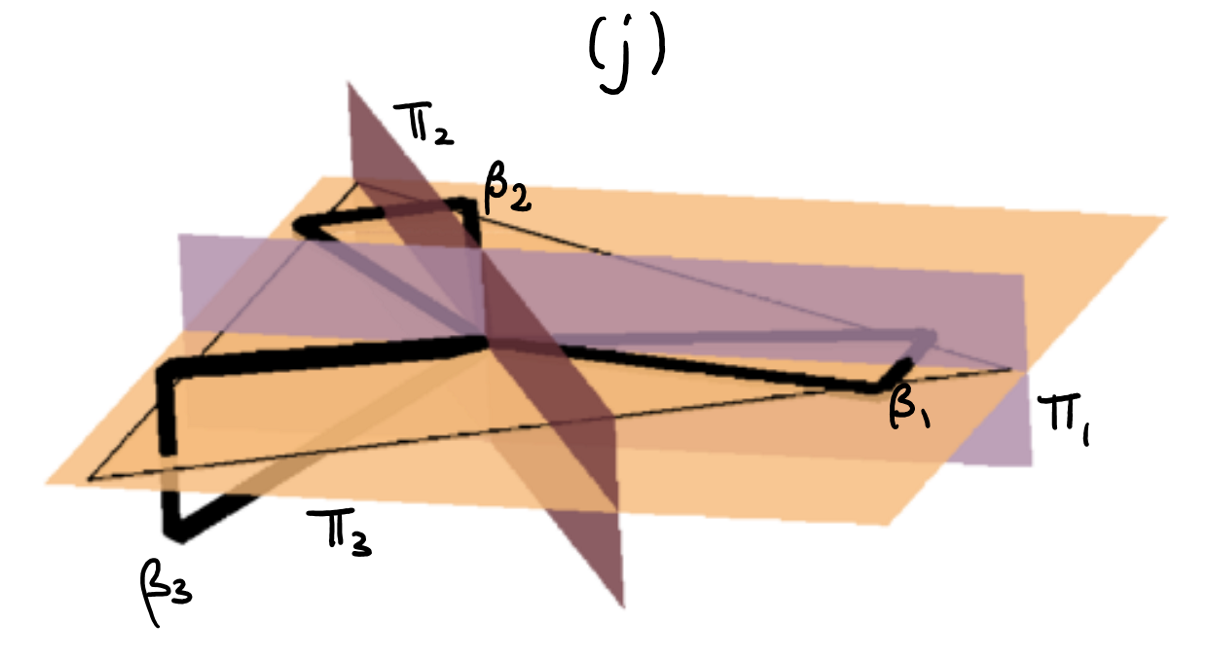}

    \caption{(d), (e) and (f): $\X_\theta$, for $\theta = \pi/4 + 0.1$, viewed from three directions. (g), (h), (i): The billiard trajectories $\beta_1$, $\beta_2$ and $\beta_3$ viewed from three directions. (j): The billiard trajectories $\beta_i$ with hyperplanes $\Pi_i = \pi(\opv{\gamma_i})$, where $\gamma_i$ is the geodesic loop in $\dbl\X_\theta$ corresponding to $\beta_i$ and $\pi : \dbl\X_\theta \to \X_\theta$ is the quotient map.}
    \label{fig:GeodesicBouquetFromSimplex_PolytopeAndBouquet}
\end{figure}

    Property (3) above implies that $\X_\theta$ contains $\{x_1,\dotsc,x_n\}$ for values of $\theta$ greater than but sufficiently close to $\pi/4$. (See \cref{fig:GeodesicBouquetFromSimplex_PolytopeAndBouquet}(d)--(f) for an illustration of $\X_\theta$ for $\theta = \pi/4 + 0.1$.) For any $\varepsilon > 0$, let $\phi = \pi/4 + \varepsilon/2$ and let $x_i^{\pm\phi}$ be the unique point in $W_i \cap \partial\HS_i^{\pm\phi} \cap \HS_i^{\mp\phi}$ such that $\angle x_i^{\pm\phi}ox_i = \varepsilon$. Now we choose $\varepsilon$ small enough such that $\{x_i^{-\phi}, x_i^\phi\} \in \HS_j^{-\phi} \cap \HS_j^\phi$ for all $1 \leq i, j \leq n$. This condition, combined with the choice of angles, implies that there is a billiard trajectory $\beta_i$ in $\X_\phi$ that travels from the origin to $x_i^{-\phi}$, and then to $x_i^\phi$ and back to the origin. (The existence of such a billiard trajectory depends on the fact that $\phi$ is strictly greater than $\pi/4$. See \cref{fig:GeodesicBouquetFromSimplex_PolytopeAndBouquet}(g)--(i) for an illustration of those billiard trajectories when $\theta = \phi = \pi/4 + 0.1$.) In fact, all of the $\beta_i$'s are congruent in the sense that they can be superimposed with one another using rigid motions. Now choose a point $p \in \dbl\X_\phi$ that maps to the origin under the quotient map $\dbl\X_\phi \to \X_\phi$. Then $\beta_i$ corresponds to a geodesic loop $\gamma_i$ in $\dbl\X_\phi$ based at $p$. Since the sum of the unit tangent vectors of $\gamma_i$ at $p$ (pointing away from $p$) is $2\cos(\varepsilon) x_i$, property (2) implies that altogether, $\gamma_1, \dotsc, \gamma_n$ forms a stationary geodesic bouquet $G$.
    
    \Cref{lem:2CollisionBilliadLoopOPV} implies that $\ker\opd{\gamma_i} = \Pi_i$ (see \cref{fig:GeodesicBouquetFromSimplex_PolytopeAndBouquet}(j)). But $\bigcap_i\Pi_i = \{0\}$ so by \cref{cor:NullvarsToIntersectKernels_Flat}, $G$ is stable.
    
    Finally, as long as we choose $\varepsilon$ to be small enough, no two tangent vectors of the loops of $G$ at the basepoint can be parallel, because $x_i$ and $x_j$ cannot be parallel when $i \neq j$. (That would violate the condition that any $n-1$ vectors among $x_1,\dotsc,x_n$ must be linearly independent.)
\end{proof}

\subsection{Proof of the main result}

In this section we will demonstrate that our main result follows from \cref{thm:StableNLoopsFlat} and the following intuitively plausible assertions about the smoothing of doubles of polytopes. These assertions will be proven in Appendices~\ref{sec:PerturbPreserveStability} and \ref{sec:SmoothingDouble}.

\begin{restatable}{proposition}{PropSmoothing}
    \label{prop:Smoothing}
    Let $\X$ be a convex $n$-polytope such that $\dbl\X$ contains a stable geodesic bouquet $G$. Then for any neighbourhood $N$ of $G$ in $\dbl[sm]\X$ that is also a compact submanifold, there exists a sequence of embeddings $\{\varphi_i : N \to M_i\}_{i = 1}^\infty$ into smooth convex hypersurfaces $M_i$ of $\R^{n+1}$ with strictly positive curvature, such that the pullback metrics $g_i$ along $\varphi_i$ from $M_i$ converge to the flat metric on $N$ in the $C^\infty$ topology.
\end{restatable}

\begin{restatable}{proposition}{PropTransferGeodesicBouquet}
    \label{prop:TransferGeodesicBouquet}
    Let $(N,g_0)$ be a compact and flat Riemannian manifold, whose interior contains a stable geodesic bouquet $G_0$ that is injective. Assume also that no two tangent vectors of the loops of $G_0$ at the basepoint are parallel. Then for any Riemannian metric $g$ on $N$ that is sufficiently close to $g_0$ in the $C^\infty$ topology, $(N,g)$ also contains a stable geodesic bouquet with the same number of loops, and in which no two tangent vectors of its loops at the basepoint are parallel.
\end{restatable}

Now we are ready to combine all of the above progress into a proof of \cref{thm:StableNLoopsPosCurv}, restated as follows.

\ThmStableNLoopsPosCurv*

\begin{proof}[Proof of \cref{thm:StableNLoopsPosCurv}]
    For each $n \geq 2$, \cref{thm:StableNLoopsFlat} gives us a compact convex $n$-polytope $\X$ whose double $\dbl\X$ contains an injective stable geodesic bouquet $G$ in which no two tangent vectors of loops of $G$ at the basepoint are parallel. It remains to apply \cref{prop:TransferGeodesicBouquet} to the sequence of embeddings obtained from \cref{prop:Smoothing}.
\end{proof}

\appendix

\section{A Linear-algebraic Lemma}
\label{sec:LinearAlgebraLemma}

In this section we will prove \cref{lem:GenPosHyperplanesFromSimplex}, restated as follows.

\LemmaGenPosHyperplanesFromSimplex*

\begin{proof}
    Choosing the $\Pi_i$'s generically should already work, but for concreteness we construct them explicitly. Let $\Pi_n = \vspan\{x_1,\dotsc,x_n\}$, and choose some $v \in X \setminus \Pi_n$. For $1 \leq i \leq n - 1$, let $\Pi_i = \vspan\{x_i, x_{i+1}, \dotsc, x_{i + n - 3}, v\}$, where the indices of the $x_i$ are taken cyclically modulo $n$. Thus $\dim \Pi_i = n - 1$. (\Cref{fig:GeodesicBouquetFromSimplex_SetupNotation}(c) illustrates this choice of $\Pi_i$'s when $n = 3$, the $x_i$'s are the vertices of an equilateral triangle in $\R^2 \times \{0\}$ that is centered at the origin, and $v = (0,0,1)$.)
    
    Then $\dim(\Pi_1 + \Pi_2) = \dim \vspan\{x_1,\dotsc,x_{n-1}, v\} = n$ because the $n-1$ vectors $x_1, \dotsc, x_{n-1}$ are linearly independent. And $\vspan\{x_2, \dotsc, x_{n-2}, v\} \subset \Pi_1 \cap \Pi_2$ but equality holds because both sides have dimension $n-2$, by the dimension formula. Similarly, $\dim((\Pi_1 \cap \Pi_2) + \Pi_3) = \dim\vspan\{x_2,\dotsc,x_n, v\} = n$ because the $n-1$ vectors $x_2, \dotsc, x_n$ are linearly independent. Similarly again $\vspan\{x_3, \dotsc, x_{n-2}, v\} \subset (\Pi_1 \cap \Pi_2) \cap \Pi_3$ but equality holds because both sides have dimension $n-3$, by the dimension formula. Continuing inductively, one can show that $\bigcap_{i=1}^{n-1}\Pi_i= \vspan\{v\}$. But then $\bigcap_{i=1}^n\Pi_i= \{0\}$.
\end{proof}

\section{The Persistence of Stable Geodesic Bouquets after Perturbing the Metric}
\label{sec:PerturbPreserveStability}

In this section we will prove \cref{prop:TransferGeodesicBouquet}, restated as follows.

\PropTransferGeodesicBouquet*

Let us outline the proof strategy and prove a lemma about the perturbation of Morse functions.

Consider the space of immersions modulo reparametrizations $\qimm{k}N$ that was defined in \cref{sec:Definitions}, where $G$ has $k$ loops. Every other Riemannian metric $g$ on $N$ induces a functional $\length_{g} : \qimm{k}N \to \R$ that gives the length of each immersion in $(N,g)$. By hypothesis, $\length_{g_0}$ has a local minimum at $[G] \in \qimm{k}N$. Intuitively, we would like to prove that if $g$ is sufficiently close to $g_0$, then $\length_{g}$ would also have a local minimum close to $[G]$.
    
To formalize this, we will define a smooth and compact finite-dimensional manifold $B$, as well as a family of embeddings $\varphi_{g} : B \to \qimm{k}N$ parametrized by metrics $g$ sufficiently close to $g_0$. These embeddings may be considered as ``finite-dimensional approximations of portions of $\qimm{k}N$,'' and they are modelled on similar constructions in the Morse theory of path space and geodesics \cite[Section~16]{Milnor_MorseTheory}. Specifically, the image of $\varphi_{g}$ consists of immersions formed by ``broken geodesics.''

We will show that $L_{g_0} = {\length_{g_0}} \circ \varphi_{g_0}$ has only one local minimum corresponding to $[G]$, which will be non-degenerate. Next, we will prove that $L_{g} = {\length_{g}} \circ \varphi_{g}$ will also have a unique non-degenerate local minimum $y$ which will correspond to the desired stable geodesic bouquet in $(N,g)$. That $y$ exists and is unique will be shown using the theory of \emph{stable mappings}, which are, roughly speaking, $C^\infty$ maps between manifolds that are ``equivalent up to changes in coordinates'' to all other maps that are sufficiently close in the $C^\infty$ topology. Formal definitions are available in  \cite{GolubitskyGuillemin_StableMappings} and the survey \cite{Ruas_OldNewResultStableMappings}, but we will only concern ourselves with the relevant implications, summarized in the following lemma.

\begin{lemma} 
    \label{lem:MorseFunctionPerturb}
    Let $B$ be a compact smooth manifold with a Morse function $f : B \to \R$ that has a unique critical point in the interior of $B$ that has index zero. Then every function $\tilde{f} : B \to \R$ that is sufficiently close to $f$ in the $C^\infty$ topology is also Morse and has a unique critical point in the interior of $B$ with index zero.
\end{lemma}

\begin{proof}
    Since $f$ is a proper Morse function, \cite[Chapter~III, Proposition~2.2]{GolubitskyGuillemin_StableMappings} and \cite[Theorem~4.1]{Mather_InfStabilityEquivStability} imply that $f$ is \emph{infinitesimally stable}. The main implication for us will be \cite[Theorem~2]{Mather_InfStabilityImpliesStability}, which guarantees that for all $\tilde{f} : B \to \R$ sufficiently close to $f$ in the $C^\infty$ topology, $\tilde{f} = h_1 \circ f \circ h_2$ for some diffeomorphisms $h_1 : \R \to \R$ and $h_2 : B \to B$. Moreover, as $\tilde{f}$ approaches $f$, $h_1$ and $h_2$ both approach the identity. As a result, $\tilde f$ has a unique critical point in the interior of $B$ that is a non-degenerate local minimum.
\end{proof}

Now we are ready to prove \cref{prop:TransferGeodesicBouquet}.

\begin{proof}[Proof of \cref{prop:TransferGeodesicBouquet}]
    A result of J.~Cheeger  \cite[Corollary~2.2]{Cheeger_InjectivityRadius} implies that for some constant $\rho > 0$ and some open neighbourhood $U$ of $g_0$, in the space of Riemannian metrics on $N$ with the $C^\infty$ topology, the injectivity radius $\inj(N,g)$ exceeds $\rho$ for all $g \in U$. Let us now construct $B$ as follows. Subdivide each edge of $G$ into arcs whose lengths (with respect to $g_0$) are shorter than $\rho/2$. This yields an embedded graph $G_0^+$ with a great number of vertices $x_0, \dotsc, x_m$, where $x_0$ is the basepoint of $G$. Consider a constant $\delta \in (0,\rho/10)$, and let $B(x_i)$ be the closed ball of radius $\delta$ in $(N,g)$ that is centered at $x_i$. Later on we will shrink the value of $\delta$ even further. Define $B$ to be the product $B(x_0) \times B(x_1)^\perp \times \dotsb \times B(x_m)^\perp$, where $B(x_i)^\perp$ consists of the points $y \in B(x_i)$ such that $y - x_i$ is orthogonal to $G_0$ at $x_i$. (This makes sense because $(N, g_0)$ is flat and $\norm{y - x_i} < \rho$.)
    
    For each $g \in U$ and every pair of vertices $x_i$ and $x_j$ that are adjacent in $G_0^+$, the fact that $\inj(N, g) > \rho$ implies that a unique minimizing geodesic in $(N, g)$ connects each $y_i \in B(x_i)$ and $y_j \in B(x_j)$. In this manner we may define an embedding $\varphi_g : B \to \qimm{k}N$ that sends $(y_0,\dotsc,y_m)$ to the equivalence class of the immersion $\bq_k \to N$ composed of those minimizing geodesics connecting $y_i$ and $y_j$ whenever $x_i$ and $x_j$ are adjacent. Observe that $\varphi_{g_0}(x_0,\dotsc,x_m) = [G_0]$. The stability of $G_0$ implies that $x = (x_0,\dotsc,x_m)$ is a non-degenerate critical point of $L_{g_0}$, because each $x_i$ for $i \neq 0$ is restricted to $B(x_i)^\perp$ and cannot be displaced a non-zero distance along a vector field that is tangent to $G_0$.
    
    The Morse lemma guarantees that $x$ is an isolated critical point, which means that we may decrease $\delta$ and thereby shrink $B$ until $x$ remains as the only critical point of $L_{g_0}$. Thus $L_{g_0}$ is a Morse function. We can choose $U$ to be small enough such that for all $g \in U$, $L_g$ is $C^\infty$-close enough to $L_{g_0}$ for us to apply \cref{lem:MorseFunctionPerturb}, which would imply that $L_g$ is also Morse and has a unique critical point $y$ in the interior of $B$ with index 0.
    
    For any $g \in U$, let $G \in \varphi_g(y)$ be the piecewise smooth immersion $\bq_k \to N$ that is formed by gluing together many geodesic arcs. Let us prove that it is a stationary geodesic bouquet, which requires us to verify that the sum of the outgoing unit tangent vectors at the basepoint sum to zero, and that adjacent arcs meeting away from the basepoint must form an angle of $\pi$. The former criterion must be satisfied, otherwise we could have reduced the value of $L_g(y)$ by perturbing $y_0$ in $B(x_0)$. To verify the latter criterion, observe that we may shrink $U$ to guarantee that if adjacent arcs meet at a vertex $y_i \in B(x_i)^\perp$, for $i \neq 0$, then the two arcs must lie on different sides of the ``disk'' $B(x_i)^\perp$. The two arcs must meet at angle $\pi$ at $y_i$; otherwise, we could have perturbed $y_i$ in some direction along $B(x_i)^\perp$ to reduce $L_g(y)$. Therefore $G$ is a stationary geodesic bouquet, whose geodesic loops intersect the disks $B(x_i)^\perp$ transversally at the points $y_i$.
    
    We may adapt \cite[Theorem~16.2]{Milnor_MorseTheory} to prove that the Hessian of $L_g$ at $y$, denoted by $\Hess_y L_g$, has the same index as $\Hess_{[G]} {\length_g}$, which must then be zero. Therefore $G$ is stable.
\end{proof}

\section{Smoothing the Double of a Convex Polytope}
\label{sec:SmoothingDouble}

In this section we will prove \cref{prop:Smoothing}, restated as follows.

\PropSmoothing*

There are well-known methods to approximate a given convex body by a sequence of smooth convex hypersurfaces \cite{Ghomi_SmoothingConvex,Minkowski_SmoothingConvex,BonnesenFenchel_ConvexBodies}. Nevertheless, we will implement our own version of this approximation to ensure that the Riemannian metrics of these hypersurfaces will converge in the $C^\infty$ topology. These hypersurfaces will be the level sets of smooth convex functions $
R^{n+1} \to \R$, which will guarantee that the level sets will be smooth and convex. However, in general the sectional curvature of convex hypersurfaces may vanish at some points. To guarantee strictly positive sectional curvature, we will consider the level sets of functions that satisfy a stronger notion of convexity borrowed from convex optimization, defined as follows.

\begin{definition}[Strong convexity  \cite{Nesterov_ConvexOptimization}]
    Given a constant $\kappa \geq 0$ and an open convex set $U \subset \R^n$, we say that a function $f : U \to \R$ is \emph{$\kappa$-strongly convex} if, for all $x_0, x_1 \in U$ and $\lambda \in [0,1]$ we have
    \begin{equation}
        f(\lambda x_0 + (1-\lambda) x_1) + \frac{\lambda(1-\lambda)}2 \kappa\norm{x_0 - x_1}^2 \leq \lambda f(x_0) + (1-\lambda)f(x_1).
    \end{equation}
\end{definition}

We note that $\kappa$-strong convexity implies continuity. Moreover, if $f$ is smooth, then $\kappa$-strong convexity is equivalent to $\ip{(\Hess_p f)v,v} \geq \kappa\norm{v}^2$ for all $p \in U$ and $v \in T_p\R^n$ \cite{Nesterov_ConvexOptimization}. As a result, the following lemma shows that a $\kappa$-strongly convex function (for any $\kappa > 0$) has level sets with strictly positive sectional curvature.

\begin{lemma}
    \label{lem:HessianPosCurv}
    Let $U$ be an open subset of $\R^n$ and $h : U \to \R$ be a smooth function whose Hessians are positive definite. Let $z$ be a regular value. Then $h^{-1}(z)$ is a smooth hypersurface whose sectional curvatures are positive.
\end{lemma}
\begin{proof}
    Let $M = h^{-1}(z)$. Let $p \in M$. Since $z$ is a regular value, the gradient of $h$ at $p$, denoted by $\grad_p h$, does not vanish. Let $u, v \in T_pM$ be orthonormal. The Gauss equation implies that
    \begin{equation*}
        K(u,v) = \frac1{\norm{\grad_p h}^2}(A(u,u)A(v,v) - A(u,v)^2),
    \end{equation*}
    where $A$ is the second fundamental form associated with the normal vector $\grad_p h$. But $A(u,v) = -\ip{\frac{\partial}{\partial u} \grad_p h,v} = -\ip{(\Hess_p h)u,v}$. Similar computations imply that $\norm{\grad_p h}^2K(u,v)$ is a $2 \times 2$ minor of $\Hess_p h$. However, by hypothesis, $\Hess_p h$ is positive definite, so the minor is positive.
\end{proof}

The functions whose levels sets will yield our desired convex hypersurfaces will be the convolutions of the squared-distance function $f(x) = \dist{x,\X \times \{0\}}^2$ with Gaussians. To prove the strong convexity of the convolution, it will help to study the Hessian of $f$ where it exists. In particular, given a convex $n$-polytope $\X$ and a point $x \in \X$, let $C(x) \subset \R^{n+1}$ be the set of points whose closest point in $\X \times \{0\}$ is $(x,0)$. Then $f$ will coincide with the squared-distance function from $(x,0)$ over $C(x)$. As shown in the next lemma, $C(x)$ will have the shape of an \emph{affine convex cone}, that is, the translation of some convex cone.

\begin{lemma}
    \label{lem:AffineConvexCone}
    For any convex $n$-polytope $\X$ for $n \geq 2$, $C(x)$ is an affine convex cone for all $x \in \X$. Moreover, if $x$ is a vertex of $\X$ then $C(x)$ has nonempty interior.
\end{lemma}
\begin{proof}
    Note that $C(x) = C_0(x) \times \R$, where $C_0(x)$ is the set of points in $\R^n$ whose closest point in $\X$ is $x$. By translating $\X$ through $\R^n$, we may assume that $x$ is at the origin. Hence it suffices to show that $C_0(x)$ is a convex cone with apex at $x$ which has nonempty interior when $x$ is a vertex of $\X$.
    
    For each $y \in C_0(x)$, the fact that $\X$ is convex and that the closest point in $\X$ to $y$ is $x$ (the origin) implies that the hyperplane through the origin that is orthogonal to the vector $y$ separates $y$ from the interior of $\X$. Let $H_y$ be the closed half-space that is bounded by this hyperplane and that contains $\X$. Then clearly $H_{\lambda y} = H_y$ for any $\lambda \geq 0$. In addition, for any $z \in C_0(x)$ and $\lambda \in [0,1]$, $H_{\lambda y + (1-\lambda)z} \supset H_y \cap H_z \supset \X$. Therefore the closest point to $\lambda y + (1-\lambda)z$ in $\X$ is also the origin, and $\lambda y + (1-\lambda)z \in C_0(x)$. That is, $C_0(x)$ is a convex cone.
    
    If $x$ is a vertex of $\X$, then there are at least $n$ supporting hyperplanes of $\X$ meeting at $x$, such that the outward-pointing normal vectors are linearly independent. Hence the parallelepiped spanned by those vectors has positive volume: its volume is equal to the absolute value of the nonzero determinant of the matrix whose columns are those vectors. Therefore $C_0(x)$, which contains this parallelepiped, has nonempty interior.
\end{proof}

Let us proceed to prove the strong convexity of certain convolutions $f * \Theta$, where $f$ is strongly convex over some affine convex cone, as will be the case in our situation. Let $B_r(z) \subset \R^n$ denote the open ball of radius $r$ and centered at $z$. For an affine convex cone $C \subset \R^n$ with apex $v$, define its \emph{projective inradius} as $\sup_{B_r(z) \subset C}r/\norm{z - v}$, where the supremum is taken over open balls of positive radius. 

\begin{lemma}
    \label{lem:SpreadStrongConvexity}
    Let $f : \R^n \to \R$ be a convex function that is $\kappa$-strongly convex when restricted to some affine convex cone $C$ with nonempty interior, for some $\kappa > 0$. Let $\Theta : \R^n \to \R$ be a smooth and radially symmetric Gaussian probability density function such that $\Theta(y) = \theta(\norm{y})$ is a function of $\norm{y}$. Let $v$ be the apex of $C$ and let $\rho$ be its projective inradius. Then for any $r > 0$, the convolution $f * \Theta$ is $\hat\kappa$-strongly convex when restricted to $B_r(0)$, for $\hat\kappa = \kappa\omega_n \theta(\frac{1 + r}\rho + \norm{v})$, where $\omega_n$ is the volume of a unit ball in $\R^n$.
\end{lemma}
\begin{proof}
    Choose any $x_0, x_1 \in B_r(0)$ and $\lambda \in [0,1]$. Choose some $z \in C$ such that $B_{1 + r}(z) \subset C$ and $(1 + r)/\norm{z - v} = \rho$. Note that for all $y \in B_1(z)$, we have $x_i + y \in B_{1 + r}(z)$ for $i = 0,1$. Thus, the convolution $\hat{f} = f * \Theta$ satisfies
    \begin{align}
        &\hat{f}(\lambda x_0 + (1-\lambda) x_1) \\
        ={}& \int_{\R^n} f(\lambda x_0 + (1-\lambda) x_1 + y) \Theta(-y)\,dy \\
        ={}& \int_{\R^n} f(\lambda (x_0 + y) + (1-\lambda) (x_1 + y)) \Theta(-y)\,dy \\
        ={}& \int_{B_1(z)} f(\lambda (x_0 + y) + (1-\lambda) (x_1 + y)) \Theta(-y)\,dy \\
        &+ \int_{\R^n \setminus B_1(z)} f(\lambda (x_0 + y) + (1-\lambda) (x_1 + y)) \Theta(-y)\,dy \\
        \leq{}&  \int_{B_1(z)} \Big(\lambda f(x_0 + y) + (1-\lambda) f(x_1 + y) - \frac{\lambda(1-\lambda)}2\kappa\norm{x_0  - x_1}^2\Big) \Theta(-y)\,dy \\
        &+ \int_{\R^n \setminus B_1(z)} \Big(\lambda f(x_0 + y) + (1-\lambda) f(x_1 + y)\Big) \Theta(-y)\,dy \\
        ={}& \int_{\R^n} \Big(\lambda f(x_0 + y) + (1-\lambda) f(x_1 + y)\Big) \Theta(-y)\,dy \\
        &- \frac{\lambda(1-\lambda)}2\norm{x_0 - x_1}^2\kappa \int_{B_1(z)}  \Theta(-y) \,dy \\
        ={}& \lambda\hat{f}(x_0) - (1-\lambda) \hat{f}(x_1) - \frac{\lambda(1-\lambda)}2\norm{x_0 - x_1}^2\kappa\int_{B_1(z)} \Theta(-y) \,dy,
    \end{align}
    so it remains to estimate $\kappa\int_{B_1(z)} \Theta(-y) \,dy$. However, since $\Theta(-y) = \theta(\norm{y})$ is a decreasing function of $\norm{y}$, and the supremum of the norm of points in $B_1(z)$ is $\norm{z} + 1$,
    \begin{equation}
        \kappa\int_{B_1(z)} \Theta(-y) \,dy \geq \kappa\int_{B_1(z)} \theta(\norm{z} + 1) \,dy = \kappa\omega_n \theta(\norm{z} + 1) \geq \kappa\omega_n \theta\left( \frac{1 + r}\rho + \norm{v} \right),
    \end{equation}
    where the last inequality holds because $\norm{z} \leq \norm{z - v} + \norm{v} = (1 + r)/\rho + \norm{v}$.
\end{proof}

Now we are ready to prove \cref{prop:Smoothing}.

\begin{proof}[Proof of \cref{prop:Smoothing}]
    Define the function $f : \R^{n+1} \to \R$ by $f(x) = \dist{x,\overline\X}^2$, where $\overline\X = \X \times \{0\}$. We will smooth $\dbl\X$ by considering the level sets of the convolution $f * \Theta$, where $\Theta$ is a radially symmetric Gaussian probability density function. Let $v$ be a vertex of $\overline\X$; \cref{lem:AffineConvexCone} guarantees that $C(v)$ is an affine convex cone with nonempty interior, thus it has a nonzero projective inradius. Over the interior of $C(v)$, $f$ coincides with the square of the distance to $(v,0)$, so its Hessian is twice of the identity matrix. As a result, $f$ is 2-strongly convex over the interior of $C(v)$. Let $\hat{f}$ denote the restriction of $f * \Theta$ to some fixed large ball containing $\overline\X$. Then \cref{lem:SpreadStrongConvexity} guarantees that $\hat{f}$ is $\kappa$-strongly convex for some $\kappa > 0$. By \cref{lem:HessianPosCurv}, the level sets of $\hat{f}$ have strictly positive sectional curvature.
    
    Choose the variance of $\Theta$ to be sufficiently small and choose a regular value $r^+$ in the image of $\hat{f}$ such that $M^+ = \hat{f}^{-1}(r^+)$ lies in a tubular neighbourhood of $\overline\X$. Let $h : M^+ \to \overline\X$ denote the projection of the tubular neighbourhood, but restricted to $M^+$. Let us find a smooth map $\varphi : N \to M^+$ that ``approximately lifts'' the map $N \hookrightarrow \dbl\X \xrightarrow{\pi} \X \xrightarrow{\id \times \{0\}} \overline\X$ over $h$. That is, the following diagram ``nearly commutes'':
    \begin{equation}
    \begin{tikzcd}[column sep=huge]
        & M^+ \arrow[d,"h"]
        \\
        N \arrow[r,"{x \mapsto (\pi(x),0)}"'] \arrow[ur,"\varphi",dashed] & \overline\X
    \end{tikzcd}
    \end{equation}
    We will then pull back metrics on $M^+$ over $\varphi$ to get the desired metrics on $N$.

    In some sense, we will break $N$ up into simpler pieces and define $\varphi$ over each piece. For each sufficiently small $\delta > 0$ and convex polytope $\Y$, let $\Y(\delta) = \{y \in \Y~|~ \dist{y,\partial\Y} \geq \delta\}$. For each face $F$ of $\X$, let $\PP_F^\delta \subset \X$ denote a prism based at $F(\delta)$ with height $\delta$. (That is, $\PP_F^\delta$ is isometric to $F(\delta) \times [0,\delta]$.) Given that $N$ is disjoint from the $(n-2)$-skeleton of $\dbl\X$, for sufficiently small $\delta$ we know that $N$ is contained inside the image of $(\X(\delta) \cup \bigcup_F \PP_F^\delta) \times \{0,1\}$ in $\dbl\X$, which we denote by $N'$.
    
    If the radius of the tubular neighbourhood and the variance of $\Theta$ are much smaller than $\delta$, then $h^{-1}(\pi(N') \times \{0\})$ is almost isometric to $N'$. This assertion can be verified separately over each $\X(\delta) \times \{i\}$ and each $\PP_F^\delta \times \{i\}$. Thus we can define $\varphi$ over $N'$ by mapping it to $h^{-1}(\pi(N') \times \{0\})$, and then restrict to $N$.
    
    Our conclusion, that we can choose a sequence of such embeddings $\varphi$ whose pullback metrics on $N$ converge to the flat metric in the $C^\infty$ topology, follows from the property that as the variance of $\Theta$ tends to 0, the functions $f * \Theta$ converge in the $C^\infty$ topology\footnote{To establish this convergence, one may begin by expressing $\Theta$ as the convolution of two Gaussians, $\Theta_1$ and $\Theta_2$. The associativity of convolution will then yield $f * \Theta = (f * \Theta_1) * \Theta_2$. The $C^\infty$ function $f * \Theta_1$ will bear all derivatives, while its convolution with $\Theta_2$ will take care of convergence, as the variances of all these Gaussians tend to zero.} after being restricted to some fixed compact neighbourhood of $\overline\X$.
\end{proof}

\section*{Acknowledgements}

The author would like to thank his academic advisors Alexander Nabutovsky and Regina Rotman for suggesting this research topic, and for valuable discussions. The author would also like to thank Isabel Beach for useful discussions.

\bibliographystyle{unsrt}
\bibliography{References}

\end{document}